\documentclass{birkmult}
\usepackage{amsmath,amsfonts,amssymb,amsthm,graphicx,url,xspace,hyperref,setspace}

\renewcommand{\th}{\ensuremath{^{\text{th}}}\xspace}

\newtheorem{thm}{Theorem}[section]
\newtheorem{lemma}[thm]{Lemma}
\newtheorem{prop}[thm]{Proposition}
\newtheorem{corr}[thm]{Corollary}
\newtheorem{q}[thm]{Question}
\newtheorem{conj}[thm]{Conjecture}
\theoremstyle{definition}
\newtheorem{ex}[thm]{Example}
\newtheorem{defn}[thm]{Definition}
\newtheorem{remark}[thm]{Remark}
\newcommand{\head}{{\rm head}}
\newcommand{\tail}{{\rm tail}}
\newcommand{\Z}{{\mathbb Z}}
\newcommand{\R}{{\mathbb R}}

\DeclareMathSymbol{\EE}{\mathbin}{AMSb}{"45}

\newcommand{\bdelta}{{\mathbf\delta}}
\newcommand{\one}{{\mathbf 1}}
\newcommand{\Lap}{\Delta}
\newcommand{\rLap}{\Lap'}
\newcommand{\wt}{\operatorname{wt}}
\newcommand{\outdeg}{\operatorname{outdeg}}
\newcommand{\indeg}{\operatorname{indeg}}

\begin{document}

\title[Chip-Firing and Rotor-Routing]{Chip-Firing and Rotor-Routing \\ on Directed Graphs}

\author[Holroyd]{Alexander E. Holroyd}
\address{Department of Mathematics \\ University of British Columbia 
}
\email{holroyd(at)math.ubc.ca}
\thanks{A. E. H. was funded in part by an NSERC discovery grant and
Microsoft Research.}

\author[Levine]{Lionel Levine}
\address{Department of Mathematics \\ University of California, Berkeley}
\email{levine(at)math.berkeley.edu}
\thanks{\noindent L. L. was supported by an NSF Graduate Research Fellowship}

\author[M\'{e}sz\'{a}ros]{Karola M\'{e}sz\'{a}ros}
\address{Department of Mathematics \\ Massachusetts Institute of Technology}
\email{karola(at)math.mit.edu}

\author[Peres]{Yuval Peres}
\address{Theory Group \\ Microsoft Research}
\email{peres(at)microsoft.com}

\author[Propp]{James Propp}
\address{Department of Mathematical Sciences \\ University of Massachusetts at Lowell}
\email{propp(at)cs.uml.edu}
\thanks{\noindent J. P. was supported by an NSF research grant.}

\author[Wilson]{David B. Wilson}
\address{Theory Group \\ Microsoft Research}
\email{dbwilson(at)microsoft.com}

\date{January 16, 2008 (revised May 9, 2013)}

\thanks{\noindent 2000 \textit{Mathematics Subject Classification.}  Primary: 82C20; secondary: 20K01, 05C25.}

\begin{abstract}
  We give a rigorous and self-contained survey of the abelian sandpile
  model and rotor-router model on finite directed graphs, highlighting
  the connections between them.  We present several intriguing open
  problems.
\end{abstract}

\maketitle

\section{Introduction}

The abelian sandpile and rotor-router models were discovered several
times by researchers in different communities operating independently.
The abelian sandpile model was invented by Dhar \cite{D} as a test-bed
for the concept of self-organized criticality introduced in
\cite{BTW}.  Related ideas were explored earlier by Engel \cite{E1,E2}
in the form of a pedagogical tool (the ``probabilistic abacus''), by
Spencer \cite[pp.~32--35]{S}, and by Lorenzini \cite{L89,L91} in
connection with arithmetic geometry.  The rotor-router model was first
introduced by Priezzhev \textit{et al.}~\cite{PDDK} (under the
name ``Eulerian walkers model'') in connection with self-organized
criticality.  It was rediscovered several times: by Rabani, Sinclair
and Wanka \cite{RSW} as an approach to load-balancing in
multiprocessor systems, by Propp \cite{Pr} as a way to derandomize
models such as internal diffusion-limited aggregation (IDLA)
\cite{DF,LBG}, and by Dumitriu, Tetali, and Winkler as part of their
analysis of a graph-based game \cite{DTW}.  Articles on the chip-firing
game in the mathematical literature include \cite{B1,B2,BLS,BL}.
Those on the rotor-router model include
\cite{thesis,intelligencer,HP,LP,LP2}.  Below we briefly describe the
two models, deferring the formal definitions to later sections.

The \textbf{abelian sandpile model} on a directed graph~$G$, also
called the \textbf{chip-firing game}, starts with a collection of chips
at each vertex of~$G$.  If a vertex $v$ has at least as many
chips as outgoing edges, it can \textbf{fire}, sending one chip
along each outgoing edge to a neighboring vertex.  After firing a
sequence of vertices in turn, the process stops when each vertex
with positive out-degree has fewer chips than out-going edges.
The order of firings does not affect the final configuration, a
fact we shall discuss in more detail in Section~\ref{defs}.

To define the \textbf{rotor-router model} on a directed graph~$G$, for
each vertex of $G$, fix a cyclic ordering of the outgoing edges.  To
each vertex $v$ we associate a \textbf{rotor} $\rho(v)$ chosen from
among the outgoing edges from $v$.  A chip performs a walk on $G$
according to the \textbf{rotor-router rule}: if the chip is at $v$, we
first increment the rotor $\rho(v)$ to its successor $e=(v,w)$ in the
cyclic ordering of outgoing edges from $v$, and then route the chip
along $e$ to $w$.  If the chip ever reaches a \textbf{sink}, i.e.\ a
vertex of $G$ with no outgoing edges, the chip will stop there;
otherwise, the chip continues walking forever.

A common generalization of the rotor-router and chip-firing models,
the \textbf{height arrow model}, was proposed in \cite{PDDK} and
studied in \cite{DR}.

We develop the basic theory of the abelian sandpile model in
section~\ref{defs} and define the main algebraic object associated
with it, the sandpile group of~$G$ \cite{D} (also called the
``critical group'' by some authors, e.g.\ \cite{B1,W}).  Furthermore,
we establish the basic results about \textbf{recurrent} chip
configurations, which play an important role in the theory.  In
Section~\ref{rotor} we define a notion of recurrent configurations for
the rotor-router model on directed graphs and give a characterization
for them in terms of oriented spanning trees of~$G$.  The sandpile
group acts naturally on recurrent rotor configurations, and this
action is both transitive and free.  We deduce appealing proofs of two
basic results of algebraic graph theory, namely the Matrix-Tree
Theorem \cite[5.6.8]{RS} and the enumeration of Eulerian tours in
terms of oriented spanning trees \cite[Cor.\ 5.6.7]{RS}.  We also
derive a family of bijections between the recurrent chip
configurations of~$G$ and the recurrent rotor configurations of~$G$.
Such bijections have been constructed before, for example in
\cite{BW}; however, our presentation differs significantly from the
previous ones.  Section~\ref{eulerian} establishes stronger results
for both models on Eulerian digraphs and undirected graphs.  In
Section~\ref{stacks} we present an alternative view of the
rotor-router model in terms of ``cycle-popping,'' borrowing an idea
from Wilson's work on loop-erased random walk; see~\cite{PW}.  We
conclude in Section~\ref{conc} by presenting some open questions.

\section{Chip-Firing} \label{defs}

In a finite directed graph (\textbf{digraph})~$G=(V,E)$, a directed
edge~$e \in E$ points from the vertex $\tail(e)$ to the vertex
$\head(e)$.  We allow self-loops ($\head(e)=\tail(e)$) as well as
multiple edges ($\head(e)=\head(e')$ and $\tail(e)=\tail(e')$) in~$G$.
The \textbf{out-degree}~$\outdeg(v)$ of a vertex~$v$ (also denoted
by~$d_v$) is the number of edges~$e$ with $\tail(e)=v$, and the
\textbf{in-degree} $\indeg(v)$ of~$v$ is the number of edges~$e$ with
$\head(e)=v$.  A vertex is a \textbf{sink} if its out-degree is zero.
A \textbf{global sink} is a sink $s$ such that from every other vertex
there is a directed path leading to~$s$.  Note that if there is a
global sink, then it is the unique sink.

If $G$ has the same number of edges from $v$ to~$w$ as from $w$
to~$v$ for all vertices $v \neq w$ then we call $G$
\textbf{bidirected}.  In particular, a bidirected graph is
obtained by replacing each edge of an undirected graph with a pair
of directed edges, one in each direction.

Label the vertices of~$G$ as $v_1,v_2,\dots,v_n$.  The
\textbf{adjacency matrix} $A$ of~$G$ is the $n\times n$ matrix whose
$(i,j)$-entry is the number of edges from~$v_i$ to~$v_j$, which we
denote by $a_{v_i,v_j}$ or $a_{ij}$.  The (graph) \textbf{Laplacian}
of~$G$ is the $n\times n$ matrix $\Lap = D - A$, where~$D$ is the
diagonal matrix whose $(i,i)$-entry is the out-degree of $v_i$, which we denote by $d_{v_i}$ or $d_i$.  That is,
        \[ \Lap_{ij} = \begin{cases}
  - a_{ij}       & \text{for $i \neq j$,} \\
  d_i - a_{ii} & \text{for $i = j$.} \end{cases} \]
Note that the entries in each row of~$\Lap$ sum to zero.  If the
vertex $v_i$ is a sink, then the $i$\th row of $\Lap$ is zero.

A \textbf{chip configuration} $\sigma$ on~$G$, also called a
\textbf{sandpile} on~$G$, is a vector of non-negative
integers indexed by the non-sink vertices of~$G$, where $\sigma(v)$
represents the number of chips at vertex~$v$.  A chip
configuration~$\sigma$ is \textbf{stable} if $\sigma(v) < d_v$ for
every non-sink vertex~$v$.  We call a vertex~$v$ \textbf{active} in
$\sigma$ if $v$ is not a sink and $\sigma(v) \geq d_v$.  An active
vertex~$v$ can \textbf{fire}, resulting in a new chip
configuration~$\sigma'$ obtained by moving one chip along each of the
$d_v$ edges emanating from~$v$; that is, $\sigma'(w) = \sigma(w) +
a_{vw}$ for all $w \neq v$ and $\sigma'(v) = \sigma(v) - d_v +
a_{vv}$.  We call the configuration~$\sigma'$ a \textbf{successor}
of~$\sigma$.

By performing a sequence of firings, we may eventually arrive at a
stable chip configuration, or we might continue firing forever, as
the following examples show.

\begin{ex} \label{bla1}
  Consider the complete directed graph on three vertices (without self-loops).
    Then placing three chips at a vertex gives a
  configuration that stabilizes in one move, while placing four chips
  at a vertex gives a configuration that never stabilizes (see
  Figure~\ref{stabilize}).
\begin{figure}[htbp]
  \centerline{\includegraphics[width=\textwidth]{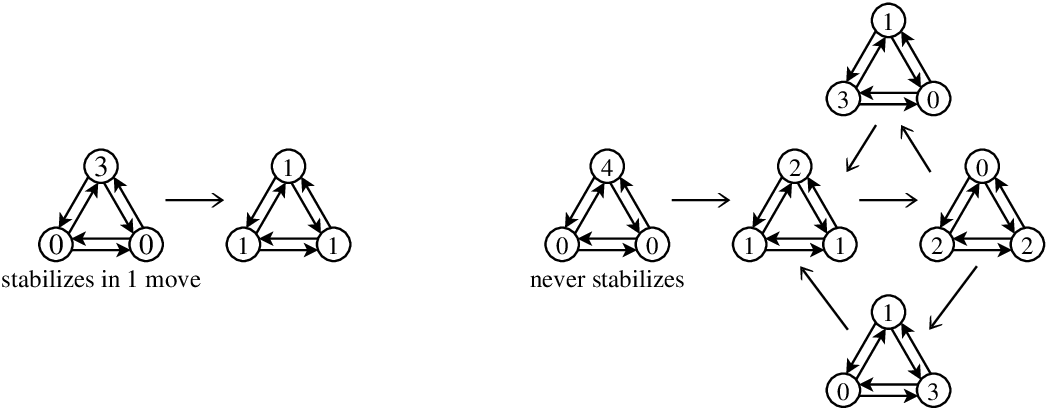}}
\caption{Some chip configurations eventually stabilize,
while others never stabilize.}
\label{stabilize}
\end{figure}
\end{ex}

\begin{figure}[htbp]
  \centerline{\includegraphics[width=\textwidth]{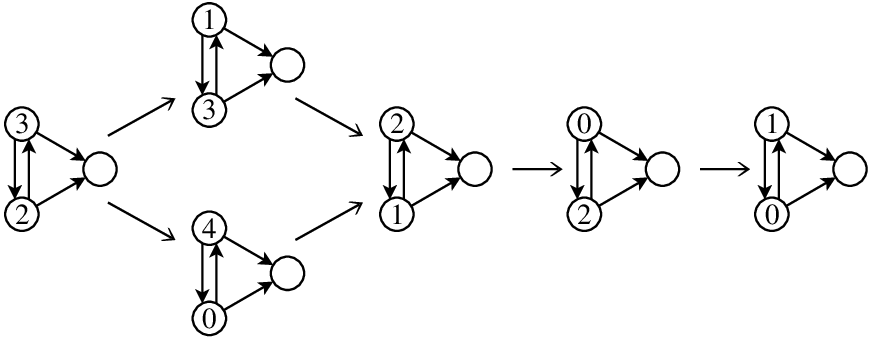}}
\caption{Commutation of the chip-firing operations.}
\label{commute}
\end{figure}

It might appear that the choice of the order in which we fire vertices
could affect the long-term behavior of the system; however, this is
not the case, as the following lemma shows (and Figure~\ref{commute}
illustrates).

\begin{lemma}[\cite{D},\cite{DF}] \label{Claim 1}
  Let $G$ be any digraph, let $\sigma_0$, $\sigma_1, \ldots, \sigma_n$
  be a sequence of chip configurations on~$G$, each of which is a
  successor of the one before, and let $\sigma'_0, \sigma'_1, \ldots,
  \sigma'_m$ be another such sequence with $\sigma'_0=\sigma_0$.
  \begin{enumerate}
  \item If $\sigma_n$ is stable, then $m \leq n$, and moreover, no
    vertex fires more times in $\sigma'_0, \ldots, \sigma'_m$ than in
    $\sigma_0, \ldots, \sigma_n$.
  \item If $\sigma_n$ and $\sigma'_m$ are both stable, then $m=n$,
    $\sigma_n=\sigma'_n$, and each vertex fires the same number of
    times in both histories.
  \end{enumerate}
\end{lemma}

\begin{proof}
  Part~2 is an immediate corollary of part~1, which we now prove.  If
  part~1 fails, then consider a counterexample with $m+n$ minimal.
  Let~$v_i$ be the vertex that fires when~$\sigma_{i-1}$
  becomes~$\sigma_i$, and~$v'_i$ be the vertex that fires
  when~$\sigma'_{i-1}$ becomes~$\sigma'_i$.  Vertex $v'_1$ must be
  fired at some stage (say the~$i$\th) in the sequence of
  configurations~$\sigma_0,\dots,\sigma_n$, since~$\sigma_n$ is
  stable; say~$v_i$ = $v'_1$.
  Then~$v_i,v_1,v_2,\dots,v_{i-1},v_{i+1},\dots,v_n$ is a permissible
  firing sequence that turns~$\sigma_0$ into~$\sigma_n$ with the same
  number of firings at each site as the unpermuted sequence.  The
  firing sequences~$v_1,\dots,v_{i-1}, v_{i+1}, \\ \dots,v_n$ and
  $v'_2,v'_3,\dots,v'_m$ then constitute a smaller counterexample to
  the lemma (with initial configuration~$\sigma'_1$), contradicting
  minimality.
\end{proof}

\begin{defn}
  Starting from a configuration~$\sigma$, Lemma~\ref{Claim 1} shows
  that there is at most one stable configuration that can be reached
  by a finite sequence of firings (and that if such a configuration
  exists, then no infinite sequence of firings is possible).  If such
  a stable configuration exists we denote it $\sigma^\circ$ and call
  it the \textbf{stabilization} of $\sigma$.
\end{defn}

Thus far, the presence or absence of sinks was irrelevant for our
claims.  For the rest of this section, we assume that the digraph~$G$
has a global sink~$s$.

\begin{lemma} \label{Claim 2}
  If digraph~$G$ has a global sink, then every chip configuration
  on~$G$ stabilizes.
\end{lemma}

\begin{proof}
  Let $N$ be the number of chips in the configuration.  Given a
  vertex $v$ of $G$, let $v_0, v_1, \dots, v_{r-1}, v_r$ be a directed
  path from $v_0 = v$ to $v_r =s$.  Every time $v_{r-1}$ fires, it
  sends a chip to the sink which remains there forever.  Thus
  $v_{r-1}$ can fire at most $N$ times.  Every time $v_{r-2}$ fires,
  it sends a chip to $v_{r-1}$, and $d_{v_{r-1}}$ such chips will
  cause $v_{r-1}$ to fire once, so $v_{r-2}$ fires at most
  $d_{v_{r-1}}N$ times.  Iterating backward along the path, we see
  that $v$ fires at most $d_{v_1}\cdots d_{v_{r-1}} N$ times.  Thus
  each vertex can fire only finitely many times, so by
  Lemma~\ref{Claim 1} the configuration stabilizes.
\end{proof}

We remark that when $G$ is connected and the sink is the only vertex
with in-degree exceeding its out-degree, the bound one gets from the
above argument on the total number of firings is far from optimal; see
Proposition~\ref{resistancebound} for a better bound.

Define the \textbf{chip addition operator}~$E_v$ as the map on
chip configurations that adds a single chip at vertex~$v$ and
then lets the system stabilize.  In symbols,
    \[ E_v \sigma = (\sigma+\one_v)^\circ \]
where $\one_v$ is the configuration consisting of a single chip at $v$.

\begin{lemma} \label{com}
  On any digraph with a global sink, the chip addition operators commute.
\end{lemma}

\begin{proof}
  Given a chip configuration~$\sigma$ and two vertices~$v$ and~$w$,
  whatever vertices are active in~$\sigma+\one_v$ are also active in
  configuration $\sigma' = \sigma+\one_v+\one_w$.  Applying to
  $\sigma'$ a sequence of firings that stabilizes $\sigma+\one_v$, we
  obtain the configuration $E_v \sigma+\one_w$.  Stabilizing this
  latter configuration yields $E_w E_v \sigma$.  Thus $E_w E_v \sigma$
  is a stabilization of $\sigma'$.  Interchanging the roles of~$v$
  and~$w$, the configuration $E_v E_w \sigma$ is also a stabilization
  of~$\sigma'$.  From Lemma~\ref{Claim 1} we conclude that $E_w E_v
  \sigma = E_v E_w \sigma$.
\end{proof}

Lemma~\ref{com} is called the \textbf{abelian property}; it
justifies Dhar's coinage ``abelian sandpile model''.  From the
above proof we also deduce the following.

\begin{corr}
  Applying a sequence of chip addition operators to~$\sigma$
  yields the same result as adding all the associated chips simultaneously
  and then stabilizing.
\end{corr}

Let $G$ be a digraph on $n$ vertices with global sink $s$.  The
\textbf{reduced Laplacian}~$\rLap$ of $G$ is obtained by deleting from
the Laplacian matrix~$\Lap$ the row and column corresponding to the
sink.  Note that firing a non-sink vertex~$v$ transforms a chip
configuration~$\sigma$ into the configuration $\sigma-\rLap_v$, where
$\rLap_v$ is the row of the reduced Laplacian corresponding to~$v$.
Since we want to view the configurations before and after firing as
equivalent, we are led to consider the group quotient $\Z^{n-1}/H$,
where $H = \Z^{n-1} \rLap$ is the integer row-span of $\rLap$.

\begin{defn} \label{groupdef}
  Let $G$ be a digraph on $n$ vertices with global sink $s$.  The
  \textbf{sandpile group} of $G$ is the group quotient
        \[ \mathcal{S}(G) = \Z^{n-1}/ \Z^{n-1}\rLap(G). \]
\end{defn}

The connection between the sandpile group and the dynamics of
sandpiles on $G$ is made explicit in Corollary~\ref{sand}.  For the
graph in Figure~\ref{zero}, the sandpile group is the cyclic group of
order~3.  The group structure of $\mathcal{S}(G)$ when $G$ is a tree
is investigated in \cite{L:tree}.

\begin{lemma} \label{size}
  The order of $\mathcal{S}(G)$ is the determinant of the reduced
  Laplacian~$\rLap(G)$.
\end{lemma}

\begin{proof}
  The order of $\mathcal{S}(G)$ equals the index of the lattice $H =
  \Z^{n-1} \rLap$ in $\Z^{n-1}$, and, recalling that the volume of a
  parallelepiped is the determinant of the matrix formed from its
  edge-vectors, we deduce that this in turn equals the determinant of
  $\rLap$.
\end{proof}

\begin{lemma} \label{chip}
  Let $G$ be a digraph with a global sink.  Every equivalence class of
  $\Z^{n-1}$ modulo~$\rLap(G)$ contains at least one stable chip
  configuration of~$G$.
\end{lemma}

\begin{proof}
  Let~$\bdelta$ be the configuration given by $\bdelta(v) = d_v$ for
  all $v$, and let $\bdelta^\circ$ be its stabilization.  Then
  $\bdelta^\circ (v) < d_v$ for all $v\neq s$, so
  $\bdelta - \bdelta^\circ$ is a positive vector equivalent to the
  zero configuration.  Given any $\alpha \in \Z^{n-1}$, let $m$
  denote the minimum of all the coordinates of $\alpha$ together
  with~0 (so that $m\leq 0$).  Then the vector
    \[ \beta = \alpha +(-m)(\bdelta - \bdelta^\circ) \]
  is nonnegative and equivalent to~$\alpha$.  Hence $\beta^\circ$ is a
  stable chip configuration in the same equivalence class as $\alpha$.
\end{proof}

\begin{ex} \label{hu}
  An equivalence class may contain more than one stable
  chip configuration.  For example, consider the complete directed
  graph on three vertices, with one of the vertices made into a sink
  by deletion of its two outgoing edges (see Figure~\ref{zero}).  It
  is easy to see that there are two stable configurations in the
  equivalence class of the identity: the configuration in which each
  of the two non-sink vertices has 0 chips and the configuration in
  which each of the two vertices has 1 chip.  It might seem natural
  that, if either of these two configurations is to be preferred as a
  representative of the identity element in the sandpile group, it
  should be the former.  However, this instinct is misleading, as we
  now explain.
\begin{figure}[htbp]
  \centerline{\includegraphics[scale=0.7]{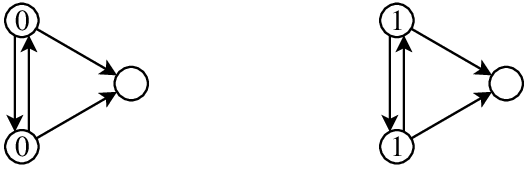}}
  \caption{Two stable chip configurations in the equivalence class
of the identity.}
\label{zero}
\end{figure}
\end{ex}

\begin{defn}
  A chip configuration~$\sigma$ is \textbf{accessible} if from any
  other chip configuration it is possible to obtain $\sigma$ by a
  combination of adding chips and selectively firing active vertices.
  A chip configuration that is both stable and accessible is called
  \textbf{recurrent}.
\end{defn}

\begin{remark}
  There are several definitions of ``recurrent'' that are used in the
  literature.  Lemma~\ref{eq} below shows that these definitions
  (including the one above) are equivalent for any digraph with a
  global sink.
\end{remark}

We will see shortly (in Lemmas~\ref{rec0} and \ref{rec1}) that each
equivalence class in $\Z^{n-1}/\Z^{n-1}\rLap$ contains a unique recurrent chip
configuration.  It is customary to represent each element of the
sandpile group by its unique recurrent element.  In Example~\ref{hu}
the all-1 configuration is accessible, but the all-0 configuration is
not.  Therefore the all-1 configuration is taken as the canonical
representative of the identity.  (However, in the context of
cluster-firing and superstabilization as described in Definition
\ref{cluster-firing} and Lemma \ref{supab}, the all-0 configuration
will be the preferred representative.)

\begin{lemma} \label{rec0}
  Let $G$ be a digraph with a global sink.  Every equivalence class of
  $\Z^{n-1}$ modulo~$\rLap(G)$ contains at least one recurrent chip
  configuration of~$G$.
\end{lemma}

\begin{proof}
  Given $\alpha \in \Z^{n-1}$, let $m$ denote the minimum of all the
  coordinates of $\alpha$ together with~0, so that $m\leq 0$.  Write
  $d_{\max}$ for the maximum out-degree of a vertex in $G$.  Then
  $\alpha$ is equivalent to the configuration
  \[ \beta = \alpha + [d_{\max}+(-m)](\bdelta-\bdelta^\circ) \]
  with $\bdelta$ as in the proof of Lemma~\ref{chip}.  Since
  $\bdelta-\bdelta^\circ$ has all entries positive, we have $\alpha \geq m(\delta-\delta^\circ)$ and hence
  $\beta \geq d_{\max}$.
  %d_{\max}(\bdelta-\bdelta^\circ) \geq d_{\max}$. 
   (Inequalities
  between two vectors or between a vector and a scalar are interpreted
  componentwise.)  In particular, $\beta$ is accessible, since any
  chip configuration can first be stabilized, so that each vertex has
  fewer than $d_{\max}$ chips, and then supplemented with extra chips
  to obtain $\beta$.  Therefore any configuration obtained from
  $\beta$ by firing is also accessible.  In particular the
  stabilization $\beta^\circ$ is thus recurrent and equivalent to
  $\alpha$.
\end{proof}

Next we will show that every equivalence class modulo $\rLap$ contains
at most one recurrent configuration, making use of the following
lemma.

\begin{lemma} \label{epsilon}
  Let $\epsilon = (2\delta)-(2\delta)^\circ$, where $\delta$ is given
  by $\delta(v) = d_v$ as before.  If $\sigma$ is recurrent, then
  $(\sigma+\epsilon)^\circ = \sigma$.
\end{lemma}

\begin{proof}
  If $\sigma$ is recurrent then it is accessible, so it can be reached
  from $\delta$ by adding some (non-negative) configuration~$\zeta$
  and selectively firing.  But since $\sigma$ is also stable this
  implies that $(\zeta+\delta)^\circ = \sigma$.  Consider the
  configuration
  \[ \gamma = (\zeta + \delta)+\epsilon = 2\delta +\zeta +\delta- (2\delta)^\circ. \]
  Since $\epsilon\geq 0$, we may start from $\gamma$ and fire a
  sequence of vertices that stabilizes $\zeta+\delta$, to obtain the
  configuration $\sigma+\epsilon$.  On the other hand, since
  $\delta-(2\delta)^\circ\geq 0$ we may start from $\gamma$ and fire a
  sequence of vertices that stabilizes $2\delta$, to obtain the
  configuration $(2 \delta)^\circ + \zeta + \delta - (2 \delta)^\circ
  = \zeta + \delta$, which in turn stabilizes to~$\sigma$.  By
  Lemma~\ref{Claim 1} it follows that $(\sigma+\epsilon)^\circ = \sigma$.
\end{proof}

\begin{lemma} \label{rec1}
  Let $G$ be a digraph with a global sink.  Every equivalence class of
  $\Z^{n-1}$ modulo~$\rLap(G)$ contains at most one recurrent chip
  configuration of~$G$.
\end{lemma}

\begin{proof}
  Let $\sigma_1$ and $\sigma_2$ be recurrent and equivalent mod
  $\rLap$.  Label the non-sink vertices $v_1, \ldots, v_{n-1}$.  Then
  $\sigma_1=\sigma_2+\sum_{i\in J} c_i \rLap_i$, where the~$c_i$ are
  nonzero constants, $\rLap_i$ is the row of the reduced
  Laplacian~$\rLap$ corresponding to $v_i$, and the index~$i$ runs
  over some subset~$J$ of the integers $1,\ldots,n-1$.  Write
  $J=J_-\cup J_+$, where $J_-=\{i:c_i<0\}$ and $J_+=\{i:c_i>0\}$, and
  let
  \[ \sigma=\sigma_1+ \sum_{i\in J_-} (-c_i) \rLap_i
     =\sigma_2+ \sum_{i\in J_+} c_i \rLap_i. \]
  Let $\epsilon$ denote the everywhere-positive chip
  configuration defined in Lemma~\ref{epsilon}.  Take $k$ large enough
  so that $\sigma' = \sigma + k\epsilon$ satisfies
  $\sigma'(v_i) \geq |c_i|d_{v_i}$ for all $i$.  Starting from
  $\sigma'$, we may fire each vertex $v_i$ for $i \in J_-$ a total of
  $-c_i$ times, and each of the intermediate configurations is a valid
  chip configuration because all the entries are nonnegative.
  The resulting configuration $\sigma_1+k\epsilon$ then
  stabilizes to $\sigma_1$ by Lemma~\ref{epsilon}.  Likewise,
  starting from $\sigma'$ we may fire each vertex $v_i$ for $i \in J_+$
  a total of $c_i$ times to obtain $\sigma_2+k\epsilon$, which
  stabilizes to $\sigma_2$.  By Lemma~\ref{Claim 1} it follows that
  $\sigma_1=\sigma_2$.
\end{proof}

\begin{corr} \label{sand}
  Let $G$ be a digraph with a global sink.  The set of all recurrent
  chip configurations on $G$ is an abelian group under the operation
  $(\sigma,\sigma') \mapsto (\sigma+\sigma')^\circ$, and it is
  isomorphic via the inclusion map to the sandpile
  group~$\mathcal{S}(G)$.
\end{corr}
\begin{proof}
  Immediate from Lemmas~\ref{rec0} and \ref{rec1}.
\end{proof}

In view of this isomorphism, we will henceforth use the term
``sandpile group'' to refer to the group of recurrent configurations.

It is of interest to consider the \textbf{identity element} $I$ of the
sandpile group, i.e.\ the unique recurrent configuration equivalent to
the all-0 configuration.  Here is one method to compute~$I$.
Let~$\sigma$ be the configuration $2\delta - 2$.  (Arithmetic
combinations of vectors and scalars are to interpreted componentwise.)
Since $\sigma^\circ \leq \delta - 1$ we have $\sigma - \sigma^\circ
\geq \delta - 1$, so $\sigma - \sigma^\circ$ is accessible.  Since
$\sigma - \sigma^\circ$ is equivalent to 0, the identity element is
given by $I = (\sigma-\sigma^\circ)^\circ$.

\begin{figure}[htbp]
 \centerline{\includegraphics[width=.48\textwidth]{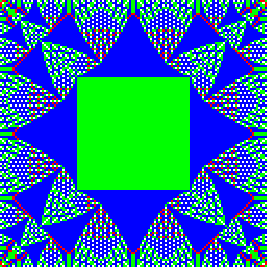}\hfil
             \includegraphics[width=.48\textwidth]{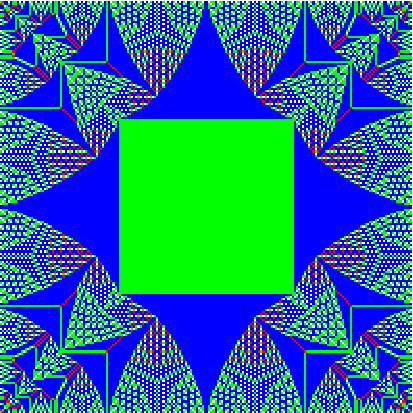} }
 \vspace*{3bp}
 \centerline{\includegraphics[width=.48\textwidth]{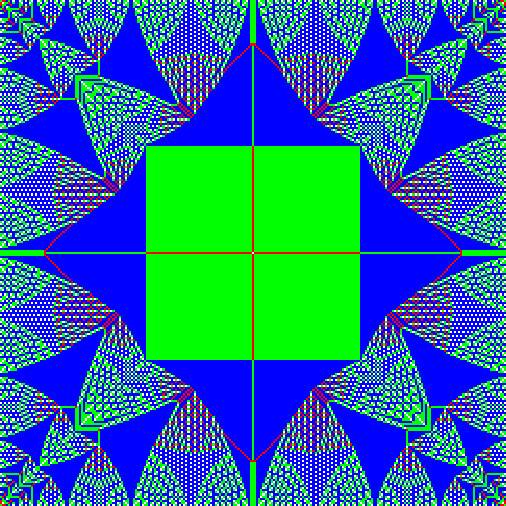}\hfil
             \includegraphics[width=.48\textwidth]{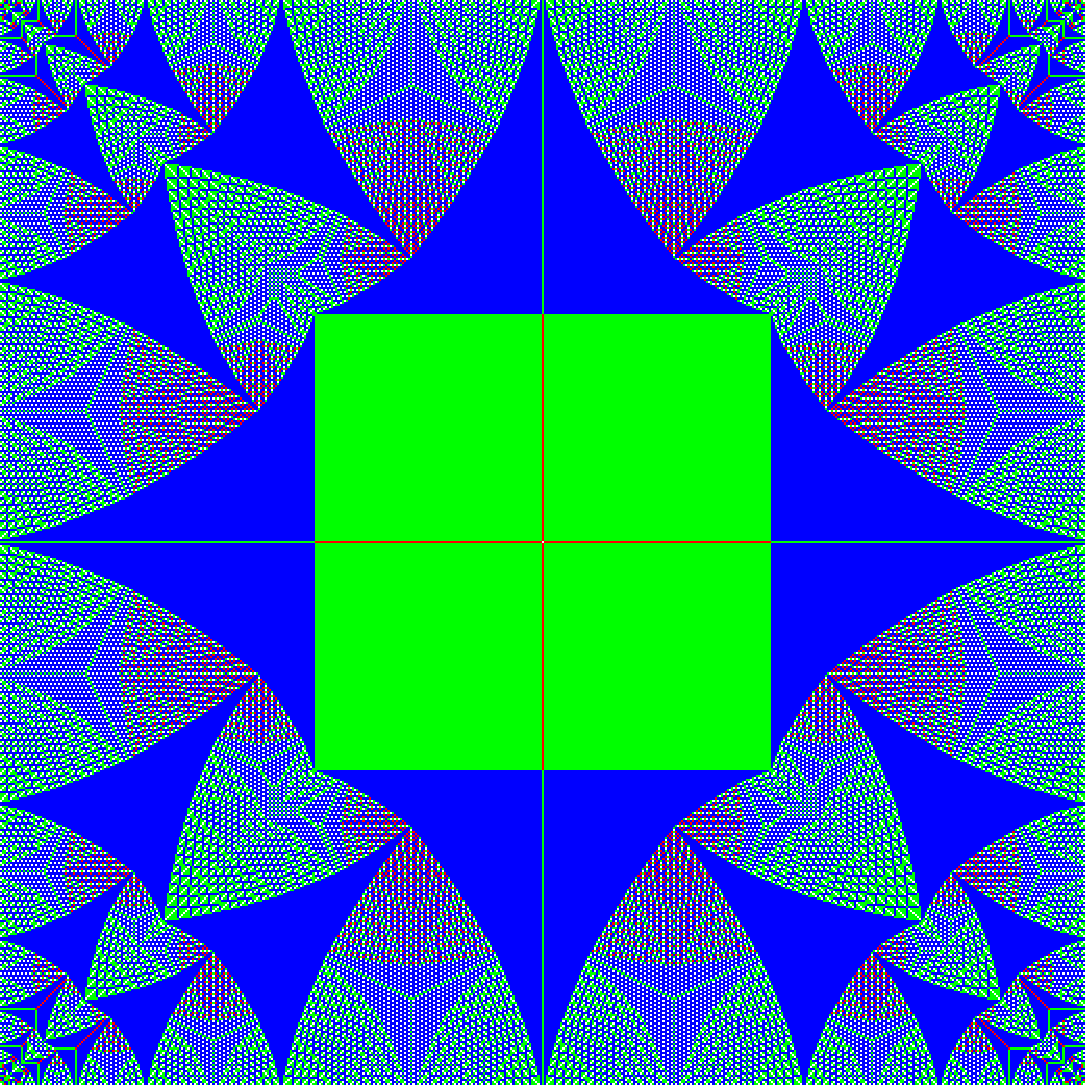} }
\caption{ The identity element of the sandpile group of the
$L\times L$ square grid for different values of $L$, namely
$L=128$ (upper left), $198$ (upper right), $243$ (lower left),
 and $521$ (lower right).
The color scheme is as follows:
orange=$0$ chips, red=$1$ chip, green=$2$ chips, and blue=$3$ chips.
 } \label{square-ident}
\end{figure}

Figure~\ref{square-ident} shows identity elements for the $L\times L$
square grid with ``wired boundary,'' for several values of~$L$.  (To
be more precise, the graph~$G$ is obtained by replacing each edge of
the undirected square grid with a pair of directed edges, and
adjoining a sink vertex~$s$ along with two edges from each of the four
corner vertices to~$s$ and one edge from each of the other boundary
vertices to~$s$.)  The identity element of this graph was studied by
Le Borgne and Rossin \cite{LBR}, but most basic properties of this
configuration, such as the existence of the large square in the
center, remain unproved.

\begin{figure}[htbp]
  \centerline{\includegraphics[width=0.4\textwidth]{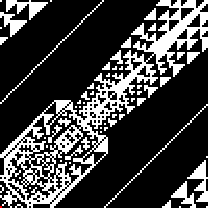}\hfil\includegraphics[width=0.4\textwidth]{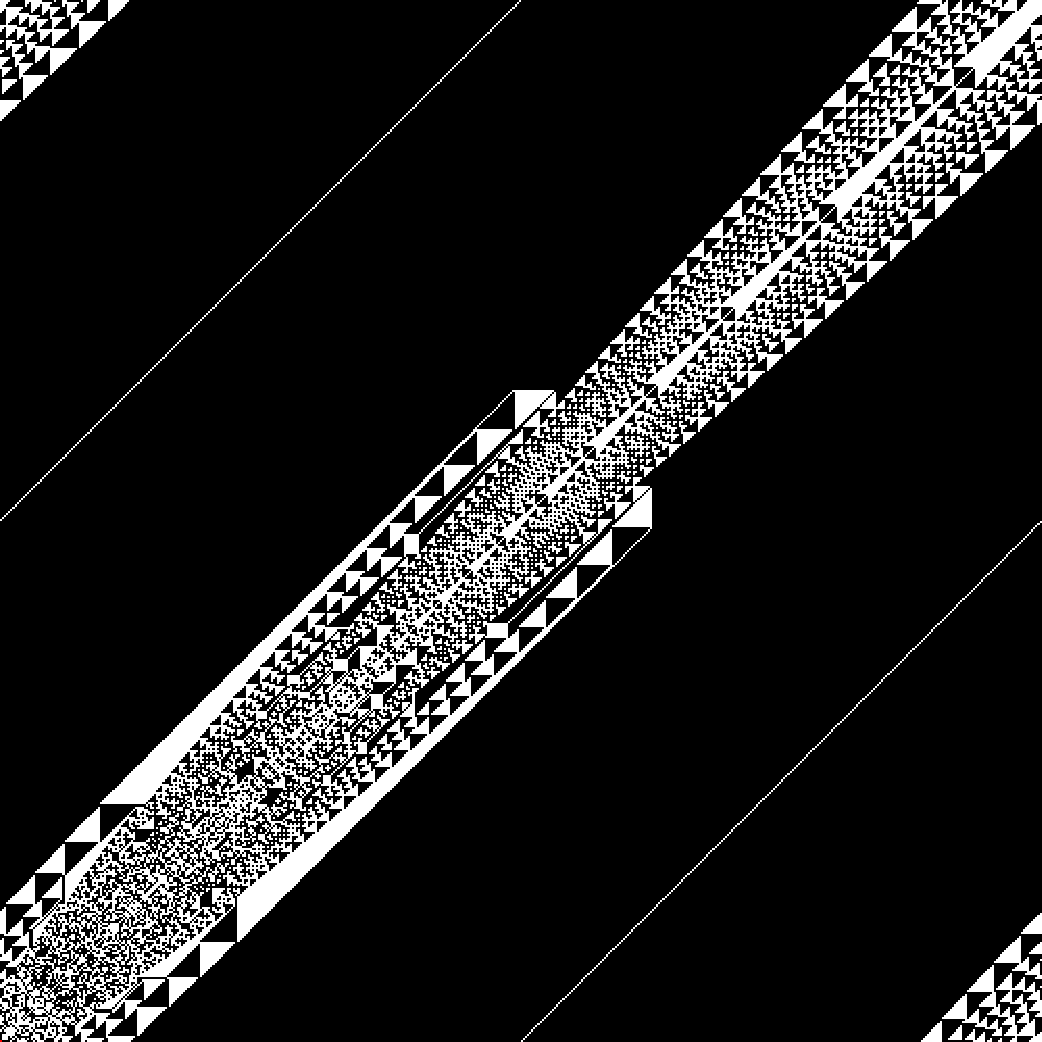}}
\caption{
The identity element of the sandpile group of the $100\times100$ directed torus (left) and the $500\times500$ directed torus (right).
The color scheme is as follows: white=$0$ chips, black=$1$ chip,
and the sink, which is at the lower-left corner, is shown in red.
}
\label{torus-ident}
\end{figure}

Figure~\ref{torus-ident} shows another example, the identity
element for the $100\times100$ directed torus.  (That is, for each
vertex $(i,j) \in \Z/100\Z \times \Z/100\Z$, there are directed
edges from $(i,j)$ to $(i+1\bmod 100,j)$ and to $(i,j+1\bmod
100)$, and we make $(0,0)$ (the lower-left vertex) into a sink
by deleting its two outgoing edges.)

\begin{figure}[htbp]
  \centerline{\includegraphics[width=0.32\textwidth]{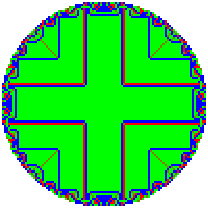}\hfil\includegraphics[width=0.32\textwidth]{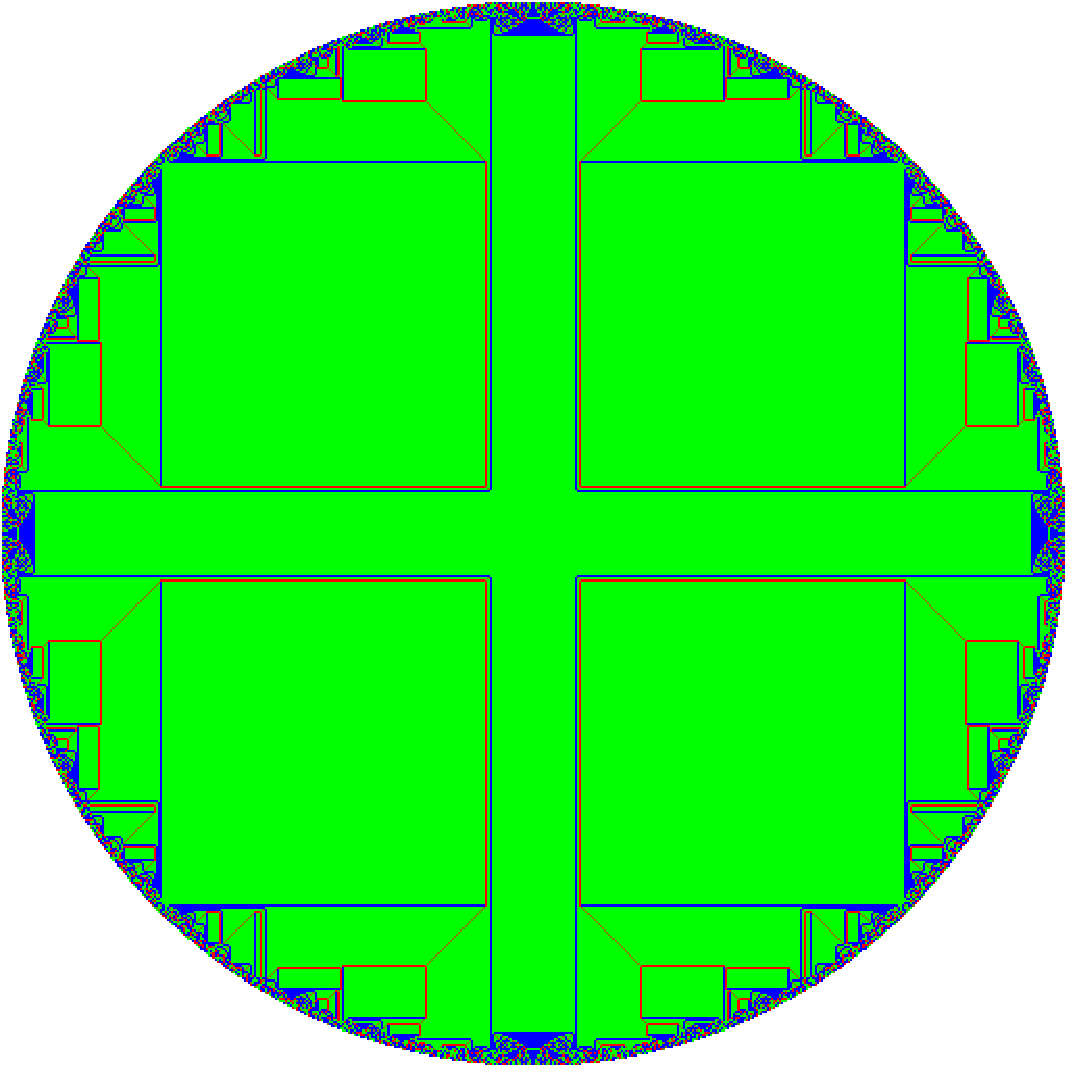}\hfil\includegraphics[width=0.32\textwidth]{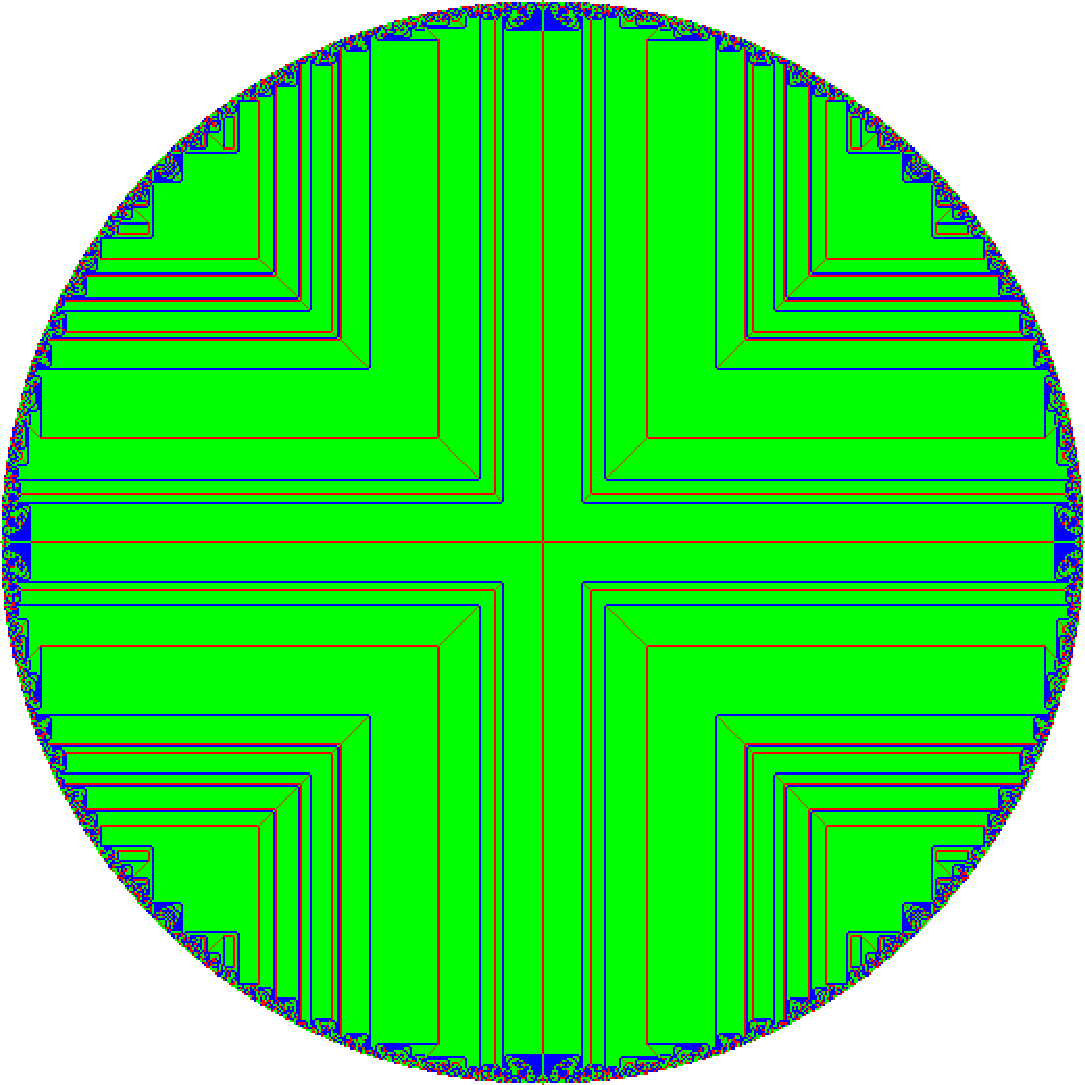}}
\caption{
  The identity element of the sandpile group of disk-shaped regions of
  diameter $100$ (left), $512$ (middle), and $521$ (right).  The color
  scheme is as follows: orange=$0$ chips, red=$1$ chip, green=$2$
  chips, and blue=$3$ chips.  }
\label{disk-ident}
\end{figure}

Figure~\ref{disk-ident} shows a third example, the identity element
for a disk-shaped region of $\Z^2$ with wired boundary.  Examples of
identity elements for graphs formed from portions of lattices other
than the square grid can be found in \cite{LKG}.

Also of interest is the \textbf{inverse} of a recurrent chip
configuration~$\sigma$ (that is, the recurrent chip
configuration~$\bar\sigma$ such that $(\sigma+\bar\sigma)^\circ=I$).  One
way to compute the inverse is via $\bar\sigma=(\zeta - \zeta^\circ -
\sigma)^\circ$, where $\zeta$ is any chip configuration satisfying
$\zeta\geq 3\delta-3$.  (Here $\zeta - \zeta^\circ - \sigma$ is
accessible, since it has at least $d_v-1$ chips at each vertex~$v$.)

Given two chip configurations $\sigma$ and $\zeta$, we say that
$\sigma$ is \textbf{reachable} from $\zeta$
(via excitation-relaxation operations) if there exists a
configuration~$\beta$ such that $\sigma = (\zeta + \beta)^\circ$.
Note that this implies that $\sigma$ is stable.  A digraph is
\textbf{strongly connected} if for any two distinct vertices $v,w$
there are directed paths from~$v$ to~$w$ and from~$w$ to~$v$.  We
write $G\smallsetminus s$ for the graph obtained from~$G$ by deleting
the vertex~$s$ along with all edges incident to~$s$.

\begin{lemma} \label{eq}
  Let $G$ be a digraph with a global sink~$s$, and let $\sigma$ be a
  chip configuration on $G$.  The following are equivalent.
  \begin{enumerate}
  \item[(1)] $\sigma$ is recurrent; that is, $\sigma$ is reachable
  from any configuration~$\zeta$.
  \item[(2)] If $\zeta$ is any configuration reachable from $\sigma$,
  then $\sigma$ is reachable from~$\zeta$.
  \item[(3)] $\sigma$ is reachable from any configuration of the form
  $E_v \sigma$, where $v$ is a non-sink vertex of~$G$.
  \item[(4)] Each strongly connected component of~$G\smallsetminus s$
  contains a vertex~$v$ such that $\sigma$ is reachable from $E_v \sigma$.
  \end{enumerate}
\end{lemma}

\begin{proof}
  Since trivially $(1) \Rightarrow (2) \Rightarrow (3) \Rightarrow (4)$,
  it suffices to show $(4)\Rightarrow(1)$.

  If (4) holds, there is a chip configuration~$\alpha$ such that
  $(\sigma+\alpha)^\circ=\sigma$ and $\alpha$ is nonzero on at least
  one vertex of each strongly connected component of~$G$.  There
  exists a positive integer~$k$ such that selective firing from $k
  \alpha$ results in a chip configuration~$\beta$ with at least one
  chip at each vertex.  Moreover, $(\sigma+\beta)^\circ = \sigma$.
  
  Now let~$\zeta$ be any chip configuration.  Since $\beta$ has at
  least one chip at each vertex, we have $\zeta \leq \sigma + \ell
  \beta$ for some integer $\ell$.  Thus we may add chips to $\zeta$ and
  then stabilize to obtain the configuration $(\sigma + \ell
  \beta)^\circ = \sigma$.  Hence $\sigma$ is recurrent.
\end{proof}

We also note that, for a digraph~$G$ with a global sink, the sandpile
group is isomorphic to the additive group of harmonic functions
modulo~1 on~$G$ that vanish on the sink \cite{Solomyak}.  A
function~$f : V(G) \to [0,1)$ is \textbf{harmonic modulo~1} if $d_v
f(v) = \sum_w a_{v,w} f(w) \mod 1$ for all vertices~$v$.  For a
sandpile configuration~$\sigma$, the associated harmonic function $f$ is the fractional part of 
the solution $\tilde{f}$ of
$$ \sum_w \rLap_{v,w} \tilde{f}(w) = \sigma(v).$$
For the graph in Figure~\ref{zero}, the three harmonic
functions are $(f(v_1),f(v_2)) = (0,0)$, $(f(v_1),f(v_2)) =
(1/3,2/3)$, and $(f(v_1),f(v_2)) = (2/3,1/3)$.

We conclude this section by pointing out a link between the sandpile
group and spanning trees.  By Lemma~\ref{size} the order of the
sandpile group of~$G$ equals the determinant of the reduced
Laplacian~$\rLap$ of~$G$.  By the \textbf{matrix-tree theorem}
\cite[5.6.8]{RS}, this determinant equals the number of
\textbf{oriented spanning trees} of~$G$ rooted at the sink (that is,
acyclic subgraphs of~$G$ in which every non-sink vertex has
out-degree~1).  Various bijections have been given for this
correspondence; see, for example, \cite{BW}.  In Section~\ref{rotor}
we will use the rotor-router model to describe a particularly natural
bijection, and deduce the matrix-tree theorem as a corollary.

\section{Rotor-Routing} \label{rotor}
 
Chip-firing is a way of routing chips through a directed graph $G$ in
such a fashion that the chips emitted by any vertex $v$ travel in
equal numbers along each of the outgoing edges.  In order to ensure
this equality, however, chips must wait at a vertex $v$ until
sufficiently many additional chips have arrived to render $v$ active.
Rotor-routing is an alternative approach to distributing chips through
$G$ which dispenses with this waiting step.  Since we cannot ensure
exact equality without waiting, we settle for the condition that the
chips emitted by any vertex $v$ travel in \textit{nearly\/} equal
numbers along each of the edges emanating from $v$.  We ensure that
this near-equality holds by using a rotor mechanism to decide where
each successive chip emitted from a vertex $v$ should be routed.

Given a directed graph~$G$, fix for each vertex $v$ a cyclic ordering
of the edges emanating from~$v$.  For an edge $e$ with tail $v$ we
denote by $e^+$ the next edge after~$e$ in the prescribed cyclic
ordering of the edges emanating from~$v$.

\begin{figure}[htbp]
  \centerline{\includegraphics[width=\textwidth]{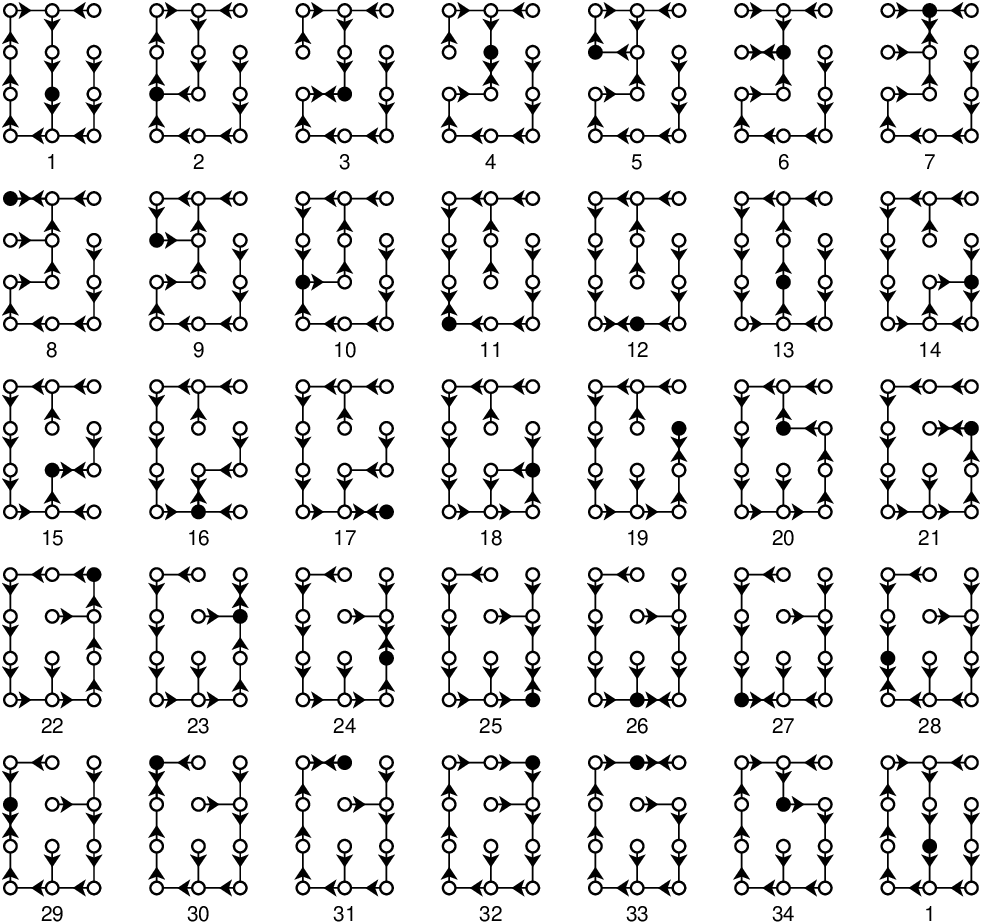}}
\caption{
  The unicycle configurations resulting from the evolution of a
  particular unicycle on the bidirected 3-by-4 rectangular grid.
  By Lemma~\ref{Claim 6}, the chip traverses each directed edge
  exactly once before the original unicycle is obtained.  Thus the
  number of distinct unicycle configurations equals the number of
  directed edges, in this case $34$.  From
  Lemma~\ref{clockwise} it follows that from any given unicycle, after some
  number of steps, the state will be the same but with the cycle's
  direction reversed.  This occurs, for example, with
  unicycles~$1$ and~$13$.}
\label{unicycle}
\end{figure}

\begin{defn}
  A \textbf{rotor configuration} is a function~$\rho$ that assigns to
  each non-sink vertex~$v$ of~$G$ an edge $\rho(v)$ emanating
  from~$v$.  If there is a chip at a non-sink vertex~$v$ of~$G$,
  \textbf{routing the chip} at $v$ (for one step) consists of updating
  the rotor configuration so that $\rho(v)$ is replaced with
  $\rho(v)^+$, and then moving the chip to the head of $\rho(v)^+$.  A
  \textbf{single-chip-and-rotor state} is a pair consisting of a
  vertex~$w$ (which represents the location of the chip) and a rotor
  configuration~$\rho$.  The \textbf{rotor-router operation} is the
  map that sends a single-chip-and-rotor state $(w,\rho)$ (where $w$
  is not a sink) to the state $(w^+,\rho^+)$ obtained by routing the
  chip at $w$ for one step.  (See Figure~\ref{unicycle} for examples
  of the rotor-router operation.)
\end{defn}
 
As we will see, there is an important link between chip-firing and
rotor-routing.  A hint at this link comes from a straightforward count
of configurations.  Recall that a stable chip configuration is a way
of assigning some number of chips between 0 and $d_v-1$ to each
non-sink vertex~$v$ of~$G$.  Thus, the number of stable configurations
is exactly $\prod_v d_v$, where the product runs over all non-sink
vertices.  This is also the number of rotor configurations on~$G$.
Other connections become apparent when one explores the appropriate
notion of recurrent states for the rotor-router model.
We will treat two cases separately: digraphs with no sink, and
digraphs with a global sink (Lemma~\ref{getthereeventually} applies to
both settings).

\begin{defn} \label{recur}
  Let $G$ be a \textbf{sink-free} digraph, i.e.\ one in which each
  vertex has at least one outgoing edge.  Starting from the state
  $(w,\rho)$, if iterating the rotor-router operation eventually leads
  back to $(w,\rho)$ we say that $(w,\rho)$ is \textbf{recurrent};
  otherwise, it is \textbf{transient}.
\end{defn}

Our first goal is to give a combinatorial characterization of the
recurrent states, Theorem~\ref{rec}.  We define a \textbf{unicycle} to
be a single-chip-and-rotor state $(w,\rho)$ for which the set of
edges~$\{\rho(v)\}$ contains a unique directed cycle, and $w$ lies on
this cycle.  (Equivalently, $\rho$ is a connected functional digraph,
and $w$ is a vertex on the unique cycle in~$\rho$.)  The following
lemma shows that the rotor-router operation takes unicycles to
unicycles.

\begin{lemma} \label{unitouni}
  Let $G$ be a sink-free digraph.  If $(w,\rho)$ is a unicycle on $G$,
  then $(w^+,\rho^+)$ is also a unicycle.
\end{lemma}

\begin{proof}
  Since $(w,\rho)$ is a unicycle, the set of edges $\{\rho(v)\}_{v
    \neq w} = \{\rho^+(v)\}_{v\neq w}$ contains no directed cycles.
  The set of edges $\{\rho^+(v)\}$ forms a subgraph of $G$ in which
  every vertex has out-degree one, so it contains a directed cycle.
  Since any such cycle must contain the edge $\rho^+(w) = \rho(w)^+$,
  this cycle is unique, and $w^+$ lies on it.
\end{proof}

\begin{lemma} \label{uniperm}
  Let $G$ be a sink-free digraph.  The rotor-router operation is a
  permutation on the set of unicycles of~$G$.
\end{lemma}

\begin{proof}
  Since the set of unicycles is finite, by Lemma~\ref{unitouni} it is
  enough to show surjectivity.  Given a unicycle $U=(w,\rho)$, let
  $U^- = (w^-,\rho^-)$ be the state obtained by moving the chip
  from~$w$ to its predecessor~$w^-$ in the unique cycle through~$w$,
  and replacing the rotor at $w^-$ with its predecessor in the cyclic
  ordering of outgoing edges from~$w^-$.  Then the rotor-router
  operation applied to $U^-$ yields $U$.  It remains to show that
  $U^-$ is a unicycle; for this it suffices to show that every
  directed cycle in~$\rho^-$ passes through~$w^-$.  Suppose that there
  is a directed cycle of rotors in~$\rho^-$ which avoids $w^-$.  Since
  $\rho^-$ agrees with $\rho$ except at~$w^-$, this same directed cycle
  occurs within $\rho$ and avoids $w^-$, a contradiction since $w^-$
  is on $\rho$'s unique cycle.
\end{proof}

\begin{corr} \label{unirec}
  Let $G$ be a sink-free digraph.  If $(w,\rho)$ is a unicycle on $G$,
  then $(w,\rho)$ is recurrent.
\end{corr}

In Lemma~\ref{recuni}, below, we show that the converse holds when $G$
is strongly connected.  We will need the following lemma, which is
analogous to Lemma~\ref{Claim 2} for the abelian sandpile.  A vertex
$w$ is \textbf{globally reachable} if for each other vertex $v$ there
is a directed path from $v$ to $w$.

\begin{lemma} \label{getthereeventually}
  Let $G$ be a digraph with a globally reachable vertex $w$.  For any
  starting vertex and rotor configuration, iterating the rotor-router
  operation a suitable number of times yields a state in which the
  chip is at $w$.
\end{lemma}

\begin{proof}
  Since $w$ is globally reachable, either $G$ is sink-free or $w$ is
  the unique sink.  Thus either we can iterate the rotor-router
  operation indefinitely, or the chip eventually visits $w$.  In the
  former case, since $G$ is finite, the chip visits some vertex $v$
  infinitely often.  But if $x$ is a vertex that is visited infinitely
  often and there is an edge from $x$ to $y$, then $y$ is also visited
  infinitely often.  Inducting along a path from $v$ to $w$, we
  conclude that the chip eventually visits $w$.
\end{proof}

\begin{lemma} \label{recuni}
  Let $G$ be a strongly connected digraph.  If $(w,\rho)$ is a
  recurrent single-chip-and-rotor state on $G$, then it is a unicycle.
\end{lemma}

\begin{proof}
  Since $G$ is strongly connected, every vertex is globally reachable.
  Hence by Lemma~\ref{getthereeventually}, if we start from any
  initial state and iterate the rotor-router rule sufficiently many
  times, the chip visits every vertex of $G$.
  
  Suppose $(w,\rho)$ is a recurrent state.  Once every vertex has been
  visited and we return to the state $(w,\rho)$, suppose the rotors at
  vertices $ v_1 \ldots, v_k$ form a directed cycle.  If~$w$ does not
  lie on this cycle, then for each $i$, the last time the chip was at
  $v_i$ it moved to $v_{i+1}$, and hence the edge from $v_i$ to
  $v_{i+1}$ was traversed more recently than the edge from $v_{i-1}$
  to $v_i$.  Carrying this argument around the cycle leads to a
  contradiction.  Thus, every directed cycle in the rotor
  configuration must pass through $w$.  But now if we start from $w$
  and follow the rotors, the first vertex we revisit must be $w$.
  Hence $(w,\rho)$ is a unicycle.
\end{proof}

Combining Corollary~\ref{unirec} and Lemma~\ref{recuni}, we have
proved the following.

\begin{thm} \label{rec}
  Let $G$ be a strongly connected digraph.  Then $(w,\rho)$ is a
  recurrent single-chip-and-rotor state on $G$ if and only if it is a
  unicycle.
\end{thm}

Next we consider the case when $G$ is a digraph with a global sink.
Note that we cannot apply the rotor-router operation to states in
which the chip is at the sink.  We call these \textbf{absorbing
  states}.  For any starting state, if we iterate the rotor-router
operation sufficiently many times, the chip must eventually arrive at
the sink by Lemma~\ref{getthereeventually}.

A \textbf{chip-and-rotor state} is a pair $\tau = (\sigma,\rho)$
consisting of a chip configuration~$\sigma$ and rotor configuration
$\rho$ on $G$.  A non-sink vertex is \textbf{active} in $\tau$ if it
has at least one chip.  If $v$ is active, then \textbf{firing} $v$
results in a new chip-and-rotor state given by replacing the rotor
$\rho(v)$ with $\rho(v)^+$ and moving a single chip from $v$ to the
head of $\rho(v)^+$ (and removing the chip if $\rho(v)^+$ is a sink).
We say that $\tau'$ is a \textbf{successor} of $\tau$ if it is
obtained from $\tau$ by firing an active vertex.  We say that $\tau$
is \textbf{stable} if no vertex can fire, i.e., all chips have moved
to a sink and disappeared.  The rotor-router operation has the
following abelian property analogous to Lemma~\ref{Claim 1}.

\begin{lemma} \label{ab}
  Let $G$ be a digraph with a global sink.  Let $\tau_0$, $\tau_1,
  \ldots, \tau_n$ be a sequence of chip-and-rotor states of $G$, each
  of which is a successor of the one before.  If $\tau_0, \tau'_1,
  \ldots, \tau'_m$ is another such sequence, and $\tau_n$ is stable,
  then $m \leq n$.  If in addition $\tau'_m$ is stable, then $m=n$ and
  $\tau_n=\tau'_n$, and for each vertex~$w$, the number of times $w$
  fires is the same for both histories.
\end{lemma}

\begin{proof}
  Let~$v_i$ and~$v'_i$ be the vertices that are fired in $\tau_{i-1}$
  and $\tau'_{i-1}$ to obtain $\tau_i$ and $\tau'_i$, respectively.
  We will show that if $\tau_n$ is stable and the sequences $v$ and
  $v'$ agree in the first $i-1$ terms for some $i\leq m$, then some
  permutation of $v$ agrees with $v'$ in the first $i$ terms.  Since
  $v'_i$ is active in $\tau_{i-1} = \tau'_{i-1}$, it must be active in
  $\tau_i, \tau_{i+1}, \dots$, until it is fired.  Since $\tau_n$ is
  stable, it follows that $v_j = v'_i$ for some $j>i$.  Let $j$ be the
  minimal such index.  Starting from~$\tau_0$, the vertices $v_1, v_2,
  \ldots, v_{i-1}, v_j, v_i, v_{i+1}, \ldots, v_{j-1}, v_{j+1},
  \ldots, v_n$ can be fired in that order, resulting in the same
  stable configuration~$\tau_n$.  Moreover, this sequence agrees with
  $v'$ in the first $i$ terms.
  
  By induction, it follows that the sequence $v'$ is an initial
  subsequence of a permutation of $v$.  In particular, $m \leq n$.  If
  $\tau'_m$ is also stable, by interchanging the roles of $\tau$ and
  $\tau'$, we obtain that $v'$ is a permutation of $v$.
\end{proof}

Given a non-sink vertex $v$ in $G$, the \textbf{chip addition
  operator}~$E_v$ is the map on rotor configurations given by adding a
chip at vertex~$v$ and iterating the rotor-router operation until the
chip moves to the sink.  By Lemma~\ref{ab} and the reasoning used in
the proof of Lemma~\ref{com}, the operators~$E_v$ commute.  This is
the \textbf{abelian property} of the rotor-router model.

\begin{figure}[htbp]
  \centerline{\includegraphics[width=\textwidth]{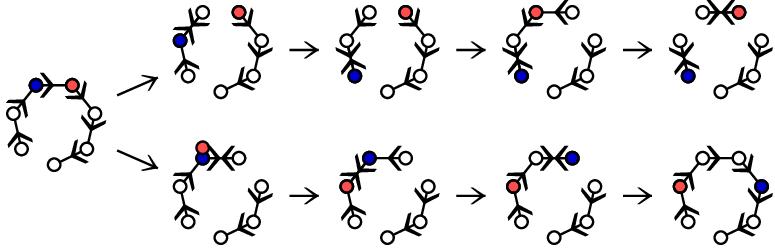}}
\caption{Failure of the abelian property for rotor-router walk
stopped after two steps.}\label{noncommute}
\end{figure}

If, rather than running the chips until they reach the sink, each chip
is run for a fixed number of steps, then the abelian property fails,
as the example in Figure~\ref{noncommute} illustrates.  (The proof of
Lemma~\ref{ab} requires that chips be indistinguishable, and it is not
possible to run each chip for a fixed number of steps without
distinguishing between them.)  Despite the failure of commutativity,
this way of routing chips has some interesting properties, similar to
the bound given in Proposition~\ref{hp}; see work of Cooper and
Spencer~\cite{CS}.

A rotor configuration~$\rho$ on $G$ is \textbf{acyclic} if the rotors
do not form any directed cycles.  If $G$ has a global sink, then
$\rho$ is acyclic if and only if the rotors form an oriented spanning
tree rooted at the sink.

\begin{lemma} \label{bij}
  Let $G$ be a digraph with a global sink, and let $v$ be a vertex
  of~$G$.  The chip addition operator~$E_v$ is a permutation on the
  set of acyclic rotor configurations on~$G$.
\end{lemma}

\begin{proof}
  We first argue that applying~$E_v$ to an acyclic rotor configuration
  yields an acyclic rotor configuration: this is proved by induction
  on the number of rotor-routing steps, where the induction hypothesis
  states that following the directed path of rotors from any vertex leads to either the sink or to the chip.
  
  Since the set of acyclic rotor configurations is finite, it suffices
  to show surjectivity.  Let $\rho$ be an acyclic rotor configuration
  on~$G$, add an edge $e$ from the sink $s$ to $v$ to form a sink-free
  digraph $G'$, and assign the rotor $\rho(s)=e$.  Then $U=(s,\rho)$
  is a unicycle on~$G'$.  Starting from $U$, we can iterate the
  \textit{inverse\/} of the rotor-router operation (which by
  Lemma~\ref{uniperm} is well-defined for unicycles) until the next
  time we reach a state with the chip at~$s$.  If we now apply the rotor-router operation once, we obtain a unicycle $U' = (v,\rho')$.
  Because $s$ belongs to the unique cycle of $U'$, deleting the edge~$e$ leaves an
  acyclic rotor configuration~$\rho'$.  Observe that running the
  rotor-router operation from $U'$ to $U$, upon ignoring the edge $e$,
  is equivalent to applying $E_v$ to $\rho'$ and obtaining $\rho$.
\end{proof}

\begin{defn}
  We next describe an action of the sandpile group on acyclic rotor
  configurations.  Given a chip configuration~$\sigma$ and a rotor
  configuration~$\rho$ on $G$, write $\sigma(\rho)$ for the rotor
  configuration obtained by adding $\sigma(v)$ chips at each vertex
  $v$ and routing them all to the sink.  By Lemma~\ref{ab} the order
  of these routings is immaterial.  Thus we may write $\sigma(\rho)$ as
  \[ \sigma(\rho) = \left( \prod_{v \in V(G) \setminus \{s\}} E_v^{\sigma(v)} \right) \rho, \]
  where the product symbol represents composition of the operators.
\end{defn}

It is trivial that $\sigma_2(\sigma_1(\rho))=(\sigma_1+\sigma_2)(\rho)$.

Since acyclic rotor configurations on $G$ can be identified with
oriented spanning trees rooted at the sink, Lemma~\ref{bij} implies
that every chip configuration~$\sigma$ acts as a permutation on the
set of oriented spanning trees of~$G$ rooted at the sink.

\begin{lemma} \label{Claim 8}
  Let $G$ be a digraph with a global sink, and let $\rho$ be an
  acyclic rotor configuration on $G$.  If the chip
  configurations~$\sigma_1$ and~$\sigma_2$ are equivalent modulo the
  reduced Laplacian $\rLap$ of $G$, then $\sigma_1(\rho) =
  \sigma_2(\rho)$.
\end{lemma}

\begin{proof}
  If $\sigma(v) \geq d_v$, and we route (for one step) $d_v$ of the
  chips at~$v$, then the rotor at~$v$ makes one full turn and one chip
  is sent along each outgoing edge from~$v$.  By Lemma~\ref{ab}, it
  follows that if $\sigma'$ is a successor to~$\sigma$ (that is,
  $\sigma'$ is obtained from~$\sigma$ by firing some active
  vertex~$v$), then $\sigma(\rho)=\sigma'(\rho)$ for any rotor
  configuration~$\rho$.  Inducting, we obtain $\sigma(\rho) =
  \sigma^\circ(\rho)$ for any rotor configuration~$\rho$.
  
  In particular, if $I$ is the recurrent chip configuration that
  represents the identity element of the sandpile group, we have
  \[ I(I(\rho)) = (I+I)(\rho) = (I+I)^\circ(\rho) = I(\rho) \]
  for any rotor configuration~$\rho$.  By Lemma~\ref{bij} the map
  $\rho \mapsto I(\rho)$ is a permutation on the set of acyclic rotor
  configurations, so it must be the identity permutation.  Now if
  $\sigma_1, \sigma_2$ are equivalent modulo $\rLap$, then
  $(\sigma_1+I)^\circ$ and $(\sigma_2+I)^\circ$ are recurrent and
  equivalent modulo $\rLap$, hence equal by Lemma~\ref{rec1}.  Since
  \[ \sigma_i(\rho) = \sigma_i(I(\rho)) = (\sigma_i+I)(\rho)
  = (\sigma_i+I)^\circ(\rho), \qquad i=1,2, \] we conclude that
  $\sigma_1(\rho)=\sigma_2(\rho)$.
 \end{proof}

It follows from Lemma~\ref{Claim 8} that the sandpile group of $G$
acts on the set of oriented spanning trees of $G$ rooted at the sink.
Our next lemma shows that this action is transitive.

\begin{lemma} \label{Claim 10}
  Let $G$ be a digraph with a global sink.  For any two acyclic rotor
  configurations~$\rho$ and~$\rho'$ on~$G$, there exists a chip
  configuration~$\sigma$ on $G$ such that $\sigma(\rho) = \rho'$.
\end{lemma}

\begin{proof}
  For a non-sink vertex $v$, let $\alpha(v)$ be the number of edges
  $e$ such that $\rho(v)<e\leq\rho'(v)$ in the cyclic ordering of
  outgoing edges from $v$.  Starting with rotor configuration~$\rho$,
  and with $\alpha(v)$ chips at each vertex~$v$, allow each chip to
  take just one step.  The resulting rotor configuration is $\rho'$;
  let $\beta$ be the resulting chip configuration, so that
  $\alpha(\rho) = \beta(\rho')$, and let $\gamma$ be the inverse of
  the corresponding element $(\beta+I)^\circ$ of the sandpile group.
  By Lemma~\ref{Claim 8} and the fact that $\beta+\gamma$ is
  equivalent to 0 modulo~$\rLap$, we have
    \[ (\alpha+\gamma)(\rho) = (\beta+\gamma)(\rho') = \rho'.  \qedhere \]
\end{proof}

Next we define recurrent rotor configurations on a digraph with a
global sink, and show they are in bijection with oriented spanning
trees.

\begin{defn} \label{recu}
  Let $G$ be a digraph with a global sink.  Given rotor
  configurations~$\rho$ and~$\rho'$ on $G$, we say that~$\rho$ is
  \textbf{reachable} from~$\rho'$ if there is a chip
  configuration~$\sigma$ such that~$\sigma(\rho')=\rho$.  We say
  that~$\rho$ is \textbf{recurrent} if it is reachable from any other
  configuration~$\rho'$.
\end{defn}

Note that in contrast to Definition~\ref{recur}, the location of the
chip plays no role in the notion of recurrent states on a digraph with
global sink.

\begin{lemma} \label{Claim 7}
  Let $G$ be a digraph with a global sink.  A rotor
  configuration~$\rho$ on~$G$ is recurrent if and only if it is
  acyclic.
\end{lemma}

\begin{proof}
  By Lemma~\ref{bij}, any configuration reachable from an acyclic
  configuration must be acyclic, so recurrent implies acyclic.
  Conversely, if $\rho$ is acyclic and $\rho'$ is any rotor
  configuration, the configuration~$\one(\rho')$ (where $\one$ denotes
  the configuration with one chip at each vertex) is acyclic, since
  the rotor at each vertex points along the edge by which a chip last
  exited.  By Lemma~\ref{Claim 10} there is a chip configuration
  $\sigma$ such that $\sigma(\one(\rho')) = \rho$, so $\rho$ is
  reachable from $\rho'$ and hence recurrent.
\end{proof}

Just as for the sandpile model, there are several equivalent
definitions of recurrence for the rotor-router model.

\begin{lemma} \label{eq1}
  Let $G$ be a digraph with a global sink $s$, and let $\rho$ be a
  rotor configuration on $G$.  The following are equivalent.
  \begin{enumerate}
  \item[(1)] $\rho$ is acyclic.
  \item[(2)] $\rho$ is recurrent; that is, $\rho$ is reachable
  from any rotor configuration~$\rho'$.
  \item[(3)] If $\rho'$ is reachable from $\rho$, then $\rho$
  is reachable from $\rho'$.
  \item[(4)] $\rho$ is reachable from any rotor configuration
  of the form $E_v \rho$, where $v$ is a vertex of $G$.
  \item[(5)] Each strongly connected component of $G\smallsetminus s$
  contains a vertex~$v$ such that $\rho$ is reachable from $E_v \rho$.
  \end{enumerate}
\end{lemma}

\begin{proof}
  By Lemma~\ref{Claim 7} we have $(1) \Rightarrow (2)$, and trivially
  $(2) \Rightarrow (3) \Rightarrow (4) \Rightarrow (5)$.
  
  If property~(5) holds, let $C_1, \ldots, C_\ell$ be the strongly
  connected components of $G\smallsetminus s$, and for each $i$, let
  $v_i \in C_i$ be such that $\rho$ is reachable from $E_{v_i} \rho$.
  Choose an integer $k$ large enough so that if we start $k$ chips at
  any $v_i$ and route them to the sink, every vertex in $C_i$ is
  visited at least once.  Let
  \[ \rho' = \left(\prod_i E_{v_i}\right)^{\!k} \rho.  \]
  Then in $\rho'$, the rotor at each vertex points along the edge by
  which a chip last exited, so $\rho'$ is acyclic.  Since $\rho$ is
  reachable from $\rho'$, by Lemma~\ref{bij} it follows that $\rho$ is
  acyclic.  Thus $(5) \Rightarrow (1)$, completing the proof.
\end{proof}

Next we show that the action of the sandpile group on the set of
oriented spanning trees of~$G$ is free.

\begin{lemma} \label{Claim 9}
  Let $G$ be a digraph with a global sink, and let $\sigma_1$ and
  $\sigma_2$ be recurrent chip configurations on~$G$.  If there is an
  acyclic rotor configuration~$\rho$ of~$G$ such that $\sigma_1(\rho) =
  \sigma_2(\rho)$, then $\sigma_1 = \sigma_2$.
\end{lemma}

\begin{proof}
  Let $\sigma = \sigma_1 + \overline{\sigma_2}$ (recall that
  $\overline{\sigma_2}$ is the inverse of $\sigma_2$.)  By
  Lemma~\ref{Claim 8}, $\sigma(\rho) =
  \overline{\sigma_2}(\sigma_1(\rho)) =
  \overline{\sigma_2}(\sigma_2(\rho)) =
  (\sigma_2+\overline{\sigma_2})(\rho)=\rho$.  Since
  $\sigma(\rho)=\rho$, after adding $\sigma$ to $\rho$, for each
  vertex $v$, the rotor at $v$ makes some integer number $c_v$ of full
  rotations.  Each full rotation results in $d_v$ chips leaving $v$,
  one along each outgoing edge.  Hence $\sigma = \sum_v c_v \rLap_v$,
  which is in the row span of the reduced Laplacian, so $\sigma$ is
  equivalent to $0$ modulo $\rLap(G)$, and hence $\sigma_1=\sigma_2$ by
  Lemma~\ref{rec1}.
\end{proof}

\begin{corr}[Matrix Tree Theorem] \label{mat}
  Let~$G$ be a digraph and $v$ a vertex of~$G$.  The number of
  oriented spanning trees of~$G$ rooted at~$v$ is equal to the
  determinant of the reduced Laplacian~$\rLap(G)$ obtained by deleting
  from $\Lap(G)$ the row and column corresponding to~$v$.
\end{corr}

\begin{proof}
  Without loss of generality we may assume the graph is loopless,
  since loops affect neither the graph Laplacian nor the number of
  spanning trees.
  
  If $v$ is not globally reachable, then there are no spanning trees
  rooted at $v$, and there is a set of vertices $S$ not containing
  $v$, such that there are no edges in $G$ from $S$ to $S^c$.  The
  rows of $\rLap(G)$ corresponding to vertices in $S$ sum to zero, so
  $\rLap(G)$ has determinant zero.
  
  If $v$ is globally reachable, delete all outgoing edges from $v$ to
  obtain a digraph~$G'$ with global sink~$v$.  Note that $G$ and $G'$
  have the same reduced Laplacian, and the same set of oriented
  spanning trees rooted at~$v$.
  
  Fix an oriented spanning tree~$\rho$ of~$G'$.  The mapping
  $\sigma\mapsto\sigma(\rho)$ from $\mathcal{S}(G')$ to the set of
  oriented spanning trees of $G'$ is a surjection by Lemma~\ref{Claim 10}
  and is one-to-one by Lemma~\ref{Claim 9}, and by
  Lemma~\ref{size}, $|\mathcal{S}(G')| = \det\rLap(G')$.
\end{proof}

Given a digraph $G$ with global sink, define its \textbf{rotor-router
  group} as the subgroup of permutations of oriented spanning trees
of~$G$ generated by the chip addition operators $E_v$.
\begin{lemma} \cite{LL}
  The rotor-router group for a digraph~$G$ with a global sink is
  isomorphic to the sandpile group $\mathcal{S}(G)$.
\end{lemma}
\begin{proof}
  The action of the sandpile group on oriented spanning trees is a
  homomorphism from the sandpile group $\mathcal{S}(G)$ onto the
  rotor-router group.  For any two distinct sandpile group elements
  $\sigma_1$ and $\sigma_2$, for any oriented spanning tree $\rho$
  (there is at least one), by Lemma~\ref{Claim 9} we have
  $\sigma_1(\rho) \neq \sigma_2(\rho)$, so the associated rotor-router
  group elements are distinct, i.e., the group homomorphism is an
  isomorphism.
\end{proof}

Since the number of recurrent chip configurations of $G$ equals
the number of oriented spanning trees, it is natural to ask for a
bijection.  Although there is no truly ``natural'' bijection,
since in general there is no canonical spanning tree to correspond
to the identity configuration, we can use the rotor-router model
to define a family of bijections.  Fix any oriented spanning
tree~$\rho$ rooted at the sink, and associate it with the identity
configuration~$I$.  For any other oriented spanning tree~$\rho'$,
by Lemma~\ref{Claim 10} there exists $\sigma \in \mathcal{S}(G)$
with $\sigma(\rho)=\rho'$; moreover, $\sigma$ is unique by
Lemma~\ref{Claim 9}.  Associate $\sigma$ with $\rho'$.  Since this
defines a surjective map from recurrent configurations to oriented
spanning trees, it must be a bijection.

\begin{remark}
  A variant of the rotor-router rule relaxes the cyclic ordering of
  edges emanating from a vertex, and merely requires one to choose
  \textit{some\/} edge emanating from the current location of the chip
  as the new rotor-setting and move the chip along this edge.  This is
  the \textbf{branching operation} introduced by Propp and studied by
  Athanasiadis~\cite{A}.  Alternatively, one can put a probability
  distribution on the edges emanating from each vertex, and stipulate
  that the new edge is to be chosen at random.  This gives the
  tree-walk introduced by Anantharam and Tsoucas~\cite{AT} in their
  proof of the Markov chain tree theorem of Leighton and
  Rivest~\cite{LR}.
\end{remark}

We conclude this section with the following result of Holroyd and
Propp \cite{HP}, which illustrates another area of application of
the rotor-router model.

\begin{prop} \label{hp}
  Let $G=(V,E)$ be a digraph and let $Y\subseteq Z$ be sets of
  vertices.  Assume that from each vertex there is a directed path
  to~$Z$.  Let $\sigma$ be a chip configuration on~$G$.  If the chips
  perform independent simple random walks on~$G$ stopped on first
  hitting~$Z$, let $H(\sigma,Y)$ be the expected number of chips that
  stop in $Y$.  If the chips perform rotor-router walks starting at
  rotor configuration~$\rho$ and stopped on first hitting~$Z$, let
  $H_\rho(\sigma,Y)$ be the number of chips that stop in~$Y$.  Then
        \begin{equation} \label{HPbound}
        | H_\rho(\sigma,Y) - H(\sigma,Y) | \leq
        \sum_{\text{\rm edges $e$}} | h(\text{\rm head of $e$})-h(\text{\rm tail of $e$}) |
        \end{equation}
      where  $h(v) := H(1_{v},Y)$.
\end{prop}
Note that the bound on the right side does not depend on~$\rho$ or~$\sigma$.

\begin{proof}
  To each edge $e=(u,v)$ with $u \in V-Z$ we assign a \textit{weight}
  \[ \wt(e) = \begin{cases} 0 &\text{if $e=\rho(u)$}, \\
    h(u)-h(v) + \wt(e^-) &\text{otherwise}.  \end{cases} \] Here $e^-$
  is the edge preceding $e$ in the cyclic ordering of edges emanating
  from~$u$.  Since $h$ is a harmonic function on~$V-Z$, the sum of $h(u)
  - h(v)$ over all edges $e=(u,v)$ emanating from~$u$ is zero, so the
  formula $\wt(e) = h(u)-h(v) + \wt(e^-)$ remains valid even when
  $e=\rho(u)$.  We assign weight $\sum_{v} \wt(\eta(v))$ to a rotor
  configuration~$\eta$, and weight~$h(v)$ to a chip located at~$v$.
  By construction, the sum of rotor and chip weights in any
  configuration is invariant under the operation of rotating the rotor
  at a chip and then routing the chip.  Initially, the sum of all chip
  weights is $H(\sigma,Y)$.  After all chips have stopped, the sum of
  the chip weights is $H_\rho(\sigma,Y)$.  Their difference is thus at
  most the change in rotor weights, which is bounded above by the sum
  in (\ref{HPbound}).
\end{proof}

Similar bounds hold even for some infinite directed graphs in
which the right side of (\ref{HPbound}) is not finite.  Thus
rotor-routing can give estimates of hitting probabilities with
very small error.  See \cite{HP} for more details.

\section{Eulerian Graphs} \label{eulerian}

A digraph $G=(V,E)$ is \textbf{Eulerian} if it is strongly connected,
and for each vertex $v \in V$ the in-degree and the out-degree of~$v$
are equal.  We call $G$ an \textbf{Eulerian digraph with sink} if it
is obtained from an Eulerian digraph by deleting all the outgoing
edges from one vertex; equivalently, $G$ has a globally reachable sink and every other
vertex has out-degree that is at least as large as its in-degree.  An
\textbf{Eulerian tour} of a digraph $G$ is a cycle in~$G$ that uses
each edge exactly once.  If $G$ has no isolated vertices, then such a tour exists if and only if $G$ is
Eulerian.  Note that for any connected undirected graph, the
corresponding bidirected graph is Eulerian.  In this section we show
some results that do not hold for general digraphs, but are true for
Eulerian ones.  We first treat the sandpile model, and then the
rotor-router model.

\begin{lemma}[Burning algorithm \cite{D}] \label{burn}
  Let~$G$ be an Eulerian digraph with sink.  A
  chip configuration~$\sigma$ is recurrent if and only if $(\sigma +
  \beta)^\circ=\sigma$, where
        \[ \beta(v) = \outdeg(v) - \indeg(v) \geq 0.  \]
  If $\sigma$ is recurrent,
  each vertex fires exactly once during the stabilization of $\sigma + \beta$.
\end{lemma}

\begin{proof}
  By the ``(4) $\Rightarrow$ (1)'' part of Lemma~\ref{eq}, if
  $(\sigma+\beta)^\circ = \sigma$, then $\sigma$ is recurrent.
  Conversely, suppose $\sigma$ is recurrent.  Label the non-sink
  vertices $v_1, \ldots, v_{n-1}$.  Since
  \begin{equation} \label{sumofrows}
  \beta = \sum_{i=1}^{n-1} \rLap_i,
  \end{equation}
  the configurations $\sigma$ and $(\sigma+\beta)^\circ$ are both
  recurrent and equivalent modulo $\rLap$.  By Lemma~\ref{rec1} it
  follows that they are equal.
  
  Let $c_i$ be the number of times vertex $v_i$ fires during the
  stabilization of $\sigma+\beta$.  Then
  \[ \sigma =
  (\sigma+\beta)^\circ = \sigma + \beta - \sum_{i=1}^{n-1} c_i
  \rLap_i. \] The rows of $\rLap$ are linearly independent, so from
  (\ref{sumofrows}) we deduce $c_i=1$ for all $i$.
\end{proof}

Informally, the burning algorithm can be described as follows: to
determine whether $\sigma$ is recurrent, first ``fire the sink'' to
obtain the configuration $\sigma+\beta$.  Then $\sigma$ is recurrent
if and only if every non-sink vertex fires in the stabilization of
$\sigma+\beta$.  In the non-Eulerian case, there is a generalization
of the burning algorithm known as the \textbf{script algorithm}, due
to Eugene Speer \cite{Speer}.

Let $H$ be an induced subgraph of $G$ not containing the sink.  We say
that $H$ is \textbf{ample} for a chip configuration~$\sigma$ on~$G$ if
there is a vertex $v$ of~$H$ that has at least as many chips as the
in-degree of $v$ in $H$.

\begin{lemma} \label{defi}
  Let $G$ be an Eulerian digraph with sink $s$.  A stable chip
  configuration~$\sigma$ on~$G$ is recurrent if and only if every
  nonempty induced subgraph of $G\smallsetminus s$ is ample for
  $\sigma$.
\end{lemma}

\begin{proof}
  If $\sigma$ is recurrent, there is a chip configuration~$\alpha$
  such that $(\bdelta+\alpha)^\circ=\sigma$, where $\bdelta(v)=d_v$.
  Each vertex of~$G$ fires at least once in the process of stabilizing
  $\bdelta+\alpha$.  Given a nonempty induced subgraph~$H$ of~$G$, let
  $v$ be the vertex of~$H$ which first finishes firing.  After $v$
  finishes firing, it must receive at least as many chips from its
  neighbors as its in-degree in $H$, so $\sigma(v)$ is at least the
  in-degree of~$v$ in~$H$.  Thus $H$ is ample for $\sigma$.
  
  Conversely, suppose that every nonempty induced subgraph of $G$ is
  ample for $\sigma$.  Let $\beta$ be the chip configuration defined
  in Lemma~\ref{burn}.  Starting from $\sigma+\beta$, fire as many
  vertices as possible under the condition that each vertex be allowed
  to fire only once.  Let $H$ be the induced subgraph on the set of
  vertices that do not fire.  Since each vertex $v$ of $H$ is unable
  to fire even after receiving one chip from each incoming edge whose
  other endpoint lies outside $H$, we have
  \[ \sigma(v) + d_v - \indeg_H(v) \leq d_v - 1.  \]
  Thus $H$ is not ample and consequently must be empty.  So every
  vertex fires once, after which we obtain the configuration
  $\sigma+\beta-\sum_{i=1}^{n-1} \rLap_i = \sigma$.  Hence
  $(\sigma+\beta)^\circ = \sigma$, which implies $\sigma$ is recurrent
  by Lemma~\ref{burn}.
\end{proof}

\begin{figure}[htbp]
  \centerline{\includegraphics[width=0.965517\textwidth]{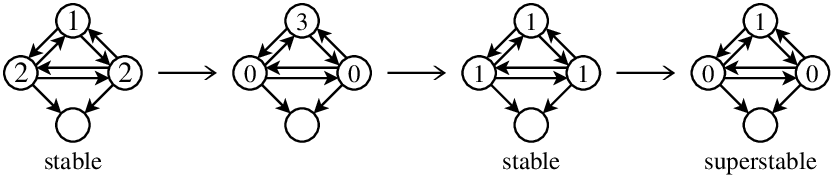}}
\caption{A sequence of cluster-firings resulting in a superstable
chip configuration.  The bottom vertex is the sink.  The clusters that fire are
first the two neighbors of the sink, next the top vertex, and finally all three
non-sink vertices.}
\label{cluster-fire}
\end{figure}

Next we define a variant of chip-firing called cluster-firing (see
Figure~\ref{cluster-fire}), and we use Lemma~\ref{defi} to
characterize the stable states for cluster-firing.  This gives rise to
a notion of ``superstable states'' which are in some sense dual to the
recurrent states.

\begin{defn} \label{cluster-firing}
  Let $G$ be a digraph with a global sink.  Let $\sigma$ be a chip
  configuration on $G$, and let $A$ be a nonempty subset of the
  non-sink vertices of $G$.  The \textbf{cluster-firing} of $A$ yields
  the configuration
  \[ \sigma' = \sigma - \sum_{i\in A} \rLap_i.  \]
  If $\sigma'$ is nonnegative we say that the cluster $A$ is
  \textbf{allowed to fire}.  We say that $\sigma$ is
  \textbf{superstable} if no cluster is allowed to fire.
\end{defn}

Note that a cluster $A$ may be allowed to fire even if no subset of
$A$ is allowed to fire.  For example, in the first configuration in
Figure~\ref{cluster-fire}, a cluster of two vertices is allowed to fire even
though the configuration is stable, so no single vertex is allowed to fire.

One could consider an even more general operation,
``multicluster-firing,'' in which different vertices can be fired
different numbers of times, so long as at the end of the firings, each
vertex has a nonnegative number of chips.  However, this
further-generalized firing operation does not yield anything new for
Eulerian digraphs, since any multicluster-firing can be expressed as a
sequence of cluster-firings: Let $m$ denote the maximal number of
times that a vertex fires in the multicluster-firing, and $C_j$ denote
the set of vertices that fire at least $j$ times in the
multicluster-firing.  Since the digraph is Eulerian, $C_m$ may be
cluster-fired, and so by induction the sets $C_m, C_{m-1}, \ldots,
C_1$ can be cluster-fired in that order.

Denote by $\delta$ the chip configuration $\delta(v)=d_v$ in which
each vertex has as many chips as outgoing edges, and by $\one$ the
configuration with a single chip at each vertex.

\begin{thm} \label{super}
  Let $G$ be an Eulerian digraph with sink.  A chip
  configuration~$\sigma$ on~$G$ is superstable if and only if $\bdelta
  - \one - \sigma$ is recurrent.
\end{thm}

\begin{proof}
  A cluster $A$ is allowed to fire if and only if for each vertex
  $v\in A$ we have
  \[ \sigma(v)-d_v+\indeg_A(v) \geq 0. \]
  This is equivalent to $d_v-1-\sigma(v)<\indeg_A(v)$, i.e., the induced
  subgraph on $A$ is not ample for $\delta-\one-\sigma$.  By
  Lemma~\ref{defi} the proof is complete.
\end{proof}

By Lemmas~\ref{rec0} and~\ref{rec1}, every equivalence class
modulo~$\rLap$ contains a unique recurrent configuration, so we obtain
the following.

\begin{corr} \label{uniquesuper}
  Let $G$ be an Eulerian digraph with sink.  Every equivalence class
  modulo~$\rLap(G)$ contains a unique superstable configuration.
\end{corr}

As a consequence, we obtain that the cluster-firing model on Eulerian
digraphs is abelian; this was proved by Paoletti \cite{P} in the
bidirected case.

\begin{corr} \label{supab}
  Let $G$ be an Eulerian digraph with sink.  Let $\sigma_0$,
  $\sigma_1, \ldots, \sigma_n$ be a sequence of chip configurations on
  $G$, each of which is obtained from the one before by a
  cluster-firing, with $\sigma_n$ superstable.  Then any sequence of
  cluster-firings that starts from $\sigma_0$ and ends in a
  superstable configuration ends in $\sigma_n$.
\end{corr}

We call the configuration~$\sigma_n$ in Corollary~\ref{supab} the
\textbf{superstabilization} of $\sigma_0$.  The following result
provides a way to compute the superstabilization.

\begin{prop}
  Let $\sigma$ be a chip configuration on an Eulerian digraph with
  sink.  The superstabilization of~$\sigma$ is given by
  \[ \sigma^* = \bdelta - \one - (\bdelta - \one - \sigma^\circ+I)^\circ \]
  where $I$ is the identity element of the sandpile group.
\end{prop}

\begin{proof}
  Since the configuration $\zeta =
  (\bdelta-\one-\sigma^\circ+I)^\circ$ is reachable from the identity
  element, it is recurrent, hence $\sigma^*=\delta-\one-\zeta$ is
  superstable by Theorem~\ref{super}.  Since $\sigma$ and $\sigma^*$
  are equivalent modulo~$\rLap$, it follows from
  Corollary~\ref{uniquesuper} that $\sigma^*$ is the
  superstabilization of~$\sigma$.
\end{proof}

\begin{figure}
\centerline{
\includegraphics[width=.48\textwidth]{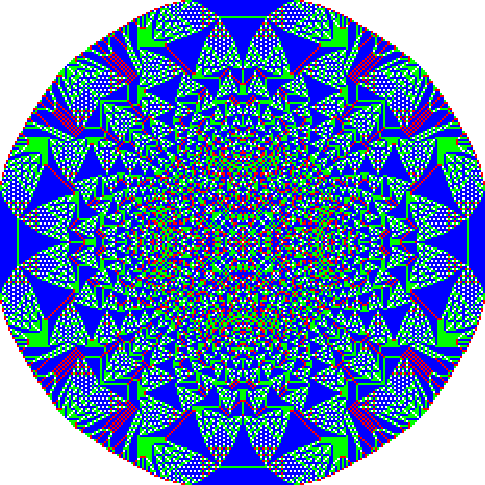}\hfil
\includegraphics[width=.48\textwidth]{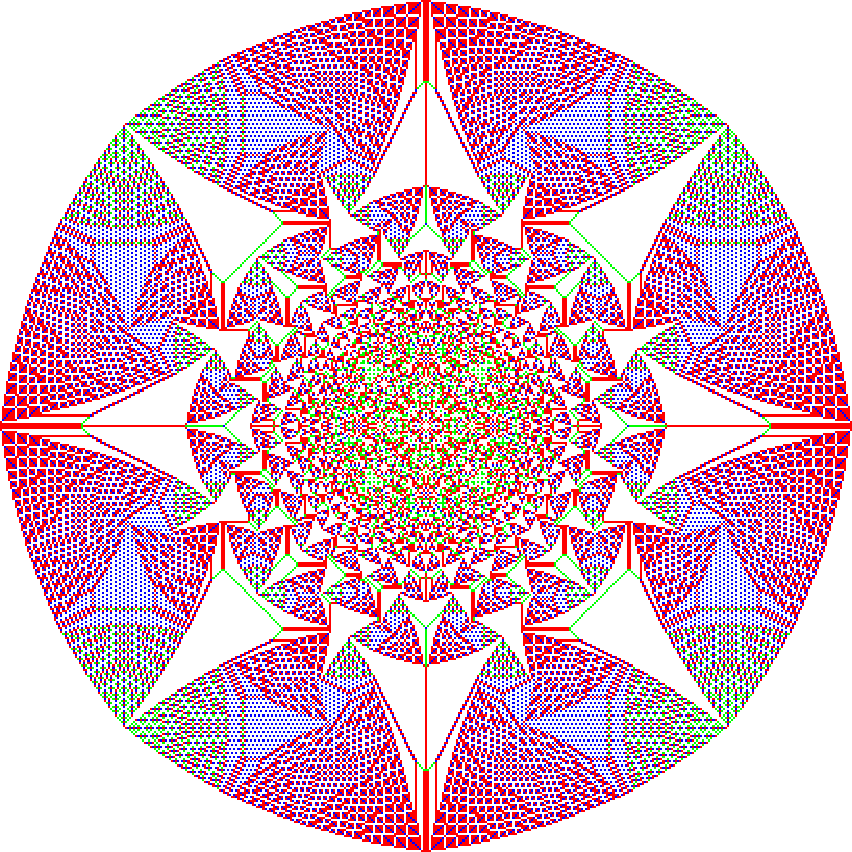}}
\caption{Stable sandpile (on left) and superstable sandpile (on right)
   of 100,000 chips, obtained by placing
  100,000 chips at the origin of the integer lattice $\Z^2$ and
  (super)stabilizing.  The color scheme is as follows: white=$0$
  chips, red=$1$ chip, green=$2$ chips, and blue=$3$ chips.  }
\end{figure}

Our final result concerning the sandpile model on Eulerian digraphs is
a theorem of Van den Heuvel \cite{VdH}; see also \cite{Tardos}.  We
give a shorter and more direct proof than that presented in
\cite{VdH}.  By an \textbf{bidirected graph with sink} $s$, we will
mean the digraph obtained from an undirected graph by first replacing
each edge by a pair of directed edges in opposite directions, and then
deleting all outgoing edges from $s$.  The \textbf{effective
  resistance} between two vertices of $G$ is an important quantity in
electrical network theory; see, e.g., \cite{DS}.  In particular, the
quantity $R_{\max}$ appearing in the proposition below is always
bounded above by the diameter of $G$, but for many graphs it is
substantially smaller than the diameter.

\begin{prop}
\label{resistancebound}
  Let~$G$ be a bidirected graph with sink, and let $\sigma$ be a chip
  configuration on~$G$.  The total number of chip moves needed to
  stabilize~$\sigma$ is at most
  \[ 2 m \,|\sigma|\, R_{\max} \]
  where~$m$ is the number of edges, $|\sigma|$ is the total number of
  chips, and $R_{\max}$ is the maximum effective resistance between any
  vertex of~$G$ and the sink.
\end{prop}

Note that firing a vertex $v$ consists of $d_v$ chip moves.

\begin{proof}
  Let $\sigma = \sigma_0, \sigma_1, \ldots, \sigma_k=\sigma^\circ$ be
  a sequence of chip configurations with $\sigma_{i+1}$ obtained from
  $\sigma_i$ by firing a single active vertex $x_i$.  Define the
  weight of~$\sigma_i$ to be
     \[ \wt(\sigma_i) = \sum_x \sigma_i(x) \wt(x), \]
  where
     \[ \wt(x) = \EE_x T_s \]
  is the expected time for a simple random walk started at~$x$ to hit the
  sink.  By conditioning on the first step $X_1$ of the walk, we compute
  \[ \Lap \wt(x) = \EE_x (\EE_{X_1} T_s - T_s) = -1, \]
  so firing the vertex~$x_i$ decreases the total weight by $d_{x_i}$.
  Thus
  \begin{equation} \label{telescopeme} \wt(\sigma_i) - \wt(\sigma_{i+1}) = d_{x_i}.  \end{equation}
  By \cite{Chandraetal}, the function $\wt$ is bounded by $2 m
  R_{\max}$.  Since the final weight $\wt(\sigma^\circ)$ is
  nonnegative, summing (\ref{telescopeme}) over $i$ we obtain that the
  total number of chip moves $N$ needed to stabilize $\sigma$ is at
  most
  \[ N = \sum_{i=0}^{k-1} d_{x_i} = \wt(\sigma) - \wt(\sigma^\circ) \leq 2m\,|\sigma|\,R_{\max}.  \qedhere \]
\end{proof}

Next we present results about the rotor-router model specific to the
Eulerian case.  An example of the next lemma is illustrated in
Figure~\ref{unicycle}.

\begin{lemma} \label{Claim 6}
  Let~$G$ be an Eulerian digraph with $m$ edges.  Let $U = (w,\rho)$
  be a unicycle in~$G$.  If we iterate the rotor-router operation~$m$
  times starting from~$U$, the chip traverses an Eulerian tour of~$G$,
  each rotor makes one full turn, and the state of the system returns
  to~$U$.
\end{lemma}

\begin{proof}
  Iterate the rotor-router operation starting from~$U$ until some
  rotor makes more than a full turn.  Let it be the rotor at
  vertex~$v$.  During this process, $v$ must emit the chip more than
  $d_v$ times.  Hence if $v \neq w$, then~$v$ must also receive the chip
  more than $d_v$ times.  Since $G$ is Eulerian, this means that some
  neighboring vertex~$u$ must send the chip to~$v$ more than once.
  However, when the chip goes from~$u$ to~$v$ for the second time, the
  rotor at~$u$ has executed more than a full turn, contradicting our
  choice of~$v$.  Thus when the rotor at~$w$ has made a full turn,
  the rotors at the other sites have made at most a full turn.
  
  We can now repeat this argument starting from the configuration
  obtained after the rotor at~$w$ has made a full turn.  In this
  way, the future history of the system is divided up into segments,
  each of length at most~$m$, where the chip is at~$w$ at the start
  of each segment.  It follows that over the course of the future
  history of the system, the chip is at~$w$ at least $d_{w}/m$ of
  the time.
  
  Since $G$ is strongly connected, we may apply this same argument to
  every state in the future history of the system, with every vertex
  of~$G$ playing the role of~$w$.  As the system evolves, the chip
  is at~$v$ at least $d_v/m$ of the time.  Since $\sum_v d_v/m = 1$,
  the chip is at~$v$ exactly $d_v/m$ of the time.  Hence, as the rotor
  at~$w$ executes a full turn, the rotors at the other sites also
  execute a full turn.  Since every rotor makes a full turn, every
  edge is traversed exactly once, so the chip traverses an Eulerian
  tour.
\end{proof}

We can use Lemma~\ref{Claim 6} to give a bijective proof of a
classical result in enumerative combinatorics relating the number of
Eulerian tours of an Eulerian graph $G$ to the number of oriented
spanning trees of $G$ (see, e.g., \cite[Cor.\ 5.6.7]{RS}).

\begin{corr}
  Let~$G=(V, E)$ be an Eulerian digraph.  Fix an edge $e \in
  E$ and let $\tail(e)=w$.  Let $\mathcal{T}(G,w)$ denote the number of
  oriented spanning trees in~$G$ rooted at~$w$, and let $\epsilon(G,e)$
  be the number of Eulerian tours in~$G$ starting with the edge~$e$.
  Then
  $$\epsilon(G,e)=\mathcal{T}(G,w) \prod_{v \in V} (d_v-1)!.$$
\end{corr}

\begin{proof}
  There are $\prod_{v \in V} (d_v-1)!$ ways to fix cyclic orderings of
  the outgoing edges from each vertex.  There are $\mathcal{T}(G,w)$
  ways to choose a unicycle $U=(w,\rho)$ with the chip at $w$ and
  the rotor $\rho(w)=e^-$, where $e^-$ is the edge preceding $e$ in
  the cyclic ordering of outgoing edges from $w$.  Given these data,
  we obtain from Lemma~\ref{Claim 6} an Eulerian tour of $G$ starting
  with the edge $e$, namely the path traversed by the rotor-router
  walk in $m$ steps.
  
  To show that this correspondence is bijective, given an Eulerian
  tour starting with the edge $e$, cyclically order the outgoing edges
  from each vertex $v$ in the order they appear in the tour.  Let
  $\rho(w)=e^-$ and for $v \neq w$ let $\rho(v)$ be the outgoing
  edge from $v$ that occurs last in the tour.  Then $U=(w,\rho)$ is
  a unicycle.
\end{proof}

The following result was first announced in \cite{PPS}, in the case of
rotor-router walk on a square lattice.

\begin{corr} \label{clockwise}
  Let~$G$ be a bidirected planar graph with the outgoing edges at
  each vertex ordered clockwise.  Let $(w,\rho)$ be a unicycle
  on~$G$ with the cycle $\mathcal{C}$ oriented clockwise.  After the
  rotor-router operation is iterated some number of times, each rotor
  internal to~$\mathcal{C}$ has performed a full rotation, each rotor
  external to~$\mathcal{C}$ has not moved, and each rotor
  on~$\mathcal{C}$ has performed a partial rotation so that
  $\mathcal{C}$ is now oriented counter-clockwise.
\end{corr}

\begin{proof}
  Let~$G'$ be the graph obtained from~$G$ by deleting all vertices and
  edges external to~$\mathcal{C}$.  Note that~$G'$, like~$G$, is
  Eulerian.  Let $\rho^-$ be the rotor configuration on $G'$ obtained from $\rho$ by reversing the rotors on $\mathcal{C}$ so that $\mathcal{C}$ is oriented counter-clockwise.
    Starting from the unicycle $U^- = (w,\rho^-)$ on $G'$ and applying the rotor-router operation~$\#\mathcal{C}$ times,
  the chip will traverse the cycle $\mathcal{C}'$, resulting in the
  state $U=(w,\rho|_{G'})$.  By Lemma~\ref{Claim 6}, further iteration
  of the rotor-router operation on $G'$ returns the system to the
  state $U^-$.  Since the outgoing edges at each vertex are ordered
  clockwise, it is straightforward to see that applying the
  rotor-router rule to~$U$ on~$G'$ and to $(w,\rho)$ on~$G$ results in
  the same evolution up until the time that state $U^-$ on $G'$ is
  reached.
\end{proof}

\begin{lemma} \label{eulerian-sink}
  Let $G$ be an Eulerian digraph, and let $G_v$ be the Eulerian
  digraph with sink obtained by deleting the outgoing edges from
  vertex $v$.  Then the abelian sandpile groups $\mathcal{S}(G_{v})$
  and $\mathcal{S}(G_{w})$ corresponding to different choices of sink
  are isomorphic.
\end{lemma}
\begin{proof}
  Recall that the sandpile group $\mathcal{S}(G_{v})$ is isomorphic to
  $\Z^{n-1}/ \Z^{n-1}\rLap(G)$; we argue that for Eulerian
  digraphs~$G$ it is also isomorphic to $\Z^{n}_0/ \Z^{n}\Lap(G)$, where $\Z^n_0$ is the set of vectors in $\Z^n$ whose coordinates sum to $0$.
  Vectors in $\Z^{n-1}$ are isomorphic to vectors in $\Z^n$ whose
  coordinates sum to~$0$, and modding out a vector in $\Z^{n-1}$ by a
  row of the reduced Laplacian~$\rLap$ corresponds to modding out the
  corresponding vector in $\Z^n$ by the corresponding row of the full
  Laplacian~$\Lap$.  For Eulerian digraphs~$G$, the last row of the
  full Laplacian~$\Lap$ is the negative of the sum of the remaining
  rows, so modding out by this extra row has no effect.
\end{proof}

We mention one other result that applies to undirected planar graphs,
due originally to Berman \cite[Prop.\ 4.1]{Berman}; see also \cite{CR}.

\begin{thm}[] \label{dual}
  If $G$ and $G^*$ are dual undirected planar graphs, then the
  sandpile groups of $G$ and $G^*$ are isomorphic.  (By
  Lemma~\ref{eulerian-sink}, the locations of the sink are
  irrelevant.)
\end{thm}

\section{Stacks and Cycle-Popping} \label{stacks}

Let $G$ be a digraph with a global sink.  In this section we describe
a more general way to define rotor-router walk on $G$, using arbitrary
stacks of rotors at each vertex in place of periodic rotor sequences.
To each non-sink vertex $v$ of $G$ we assign a bi-infinite
\textbf{stack} $\rho(v) = (\rho_k(v))_{k \in \Z}$
of outgoing edges from~$v$.  To \textbf{pop} the stack, we shift it to
obtain $(\rho_{k+1}(v))_{k\in \Z}$.  To \textbf{reverse pop} the
stack, we shift it in the other direction to obtain
$(\rho_{k-1}(v))_{k \in \Z}$.  The rotor-router walk can be defined in
terms of stacks as follows: if the chip is at vertex~$v$, pop the
stack~$\rho(v)$, and then move the chip along the edge $\rho_1(v)$.
We recover the ordinary rotor-router model in the case when each stack
$\rho(v)$ is a periodic sequence of period $d_v$ in which each
outgoing edge from $v$ appears once in each period.

The collection of stacks $\rho=(\rho(v))$, where $v$ ranges over the
non-sink vertices of $G$, is called a \textbf{stack configuration} on
$G$.  We say that $\rho$ is \textbf{infinitive} if for each edge
$e=(v,w)$, and each positive integer $K$, there exist stack elements
\[ \rho_k(v) = \rho_{k'}(v) = e \]
with $k \geq K$ and $k' \leq -K$.  This condition guarantees that
rotor-router walk eventually reaches the sink.

Given a stack configuration~$\rho$, the stack elements $\rho_0(v)$
define a rotor configuration on~$G$.  We say that $\rho$ is
\textbf{acyclic} if $\rho_0$ contains no directed cycles.  If
$\mathcal{C} = \{v_1, \ldots, v_m\}$ is a directed cycle in $\rho_0$,
define $\mathcal{C}\rho$ to be the stack configuration obtained by
reverse popping each of the stacks $\rho(v_i)$; we call this reverse
popping the cycle $\mathcal{C}$.  (If $\mathcal{C}$ is not a directed
cycle in $\rho_0$, set $\mathcal{C}\rho=\rho$.)

\begin{thm} \cite{PW} \label{cyc}
  Let $G$ be a digraph with a global sink, and let $\rho^0$ be an
  infinitive stack configuration on $G$.  There exist finitely many
  cycles $\mathcal{C}_1, \ldots, \mathcal{C}_m$ such that the stack
  configuration
  \[ \rho = \mathcal{C}_m \cdots \mathcal{C}_1 \rho^0 \]
  is acyclic.  Moreover, if $\mathcal{C}'_1, \ldots, \mathcal{C}'_n$
  is any sequence of cycles such that the stack configuration $\rho' =
  \mathcal{C}'_n \cdots \mathcal{C}'_1 \rho^0$ is acyclic, then $\rho'
  = \rho$.
\end{thm}

If $v$ is a non-sink vertex of $G$, the \textbf{chip addition
  operator} $E_v$ applied to the infinitive stack configuration~$\rho$
is the stack configuration~$\rho'$ obtained by adding a chip at $v$
and performing rotor-router walk until the chip reaches the sink.  The
next lemma shows that these operators commute with cycle-popping.

\begin{lemma} \label{poponecycle}
  Let $G$ be a digraph with a global sink, let $\rho$ be an infinitive
  stack configuration on $G$, and let $\mathcal{C}$ be a directed
  cycle in $G$.  Then
         \[ E_v(\mathcal{C} \rho) = \mathcal{C} (E_v \rho). \]
\end{lemma}

\begin{proof}
  Write $\rho' = E_v\rho$.  Let $v=v_0, v_1, \ldots, v_n = s$ be the
  path taken by a chip performing rotor-router walk from $v$ to the
  sink starting with stack configuration~$\rho$.  If this path is
  disjoint from $\mathcal{C}$, then the chip performs the same walk
  starting with stack configuration~$\mathcal{C}\rho$, and the cycle
  $\mathcal{C}$ is present in $\rho'_0$ if and only if it is present
  in $\rho_0$, so the proof is complete.
  
  Otherwise, choose $k$ minimal and $\ell$ maximal with $v_k, v_\ell
  \in\mathcal{C}$.  The rotor $\rho'_0(v_\ell)$ points to a vertex not
  in $\mathcal{C}$, so the cycle $\mathcal{C}$ is not present in
  $\rho'_0$.  Thus we must show $E_v(\mathcal{C}\rho) = \rho'$.  With
  stack configuration~$C\rho$, the chip will first travel the path
  $v_0, \ldots, v_k$, next traverse the cycle $\mathcal{C}$, and
  finally continue along the remainder of the path $v_k, \ldots, v_n$.
  Thus the stack at each vertex $w\in\mathcal{C}$ is popped one more
  time in going from $\mathcal{C}\rho$ to $E_v(\mathcal{C}\rho)$ than
  in going from $\rho$ to $E_v \rho$; the stack at each vertex $w
  \notin \mathcal{C}$ is popped the same number of times in both
  cases.
\end{proof}

The next lemma uses cycle-popping to give a constructive proof of the
injectivity of the chip addition operators $E_v$ on acyclic stack
configurations.  In the case of periodic rotor stacks, we gave a
non-constructive proof in Lemma~\ref{bij}.

\begin{lemma} \label{acyc}
  Let $G$ be a digraph with a global sink.  Given an acyclic
  infinitive stack configuration~$\rho$ on $G$ and a non-sink vertex
  $v$, there exists an acyclic infinitive stack configuration~$\rho'$
  such that $E_v \rho'=\rho$.
\end{lemma}

\begin{proof}
  Let $\rho^0$ be the stack configuration obtained from $\rho$ by
  reverse popping the stack at each of the vertices on the unique path
  in~$\rho_0$ from~$v$ to the sink.  A rotor-router walk started at
  $v$ with stack configuration~$\rho^0$ will travel directly along
  this path to the sink, so $E_v \rho^0=\rho$.  If $\rho^0$ is
  acyclic, the proof is complete.  Otherwise, by Theorem~\ref{cyc}
  there are cycles $\mathcal{C}_1, \ldots, \mathcal{C}_m$ such that
  $\rho' = \mathcal{C}_m \cdots \mathcal{C}_1 \rho^0$ is acyclic.  By
  Lemma~\ref{poponecycle}, we have
  \[ E_v \rho' = \mathcal{C}_m \cdots \mathcal{C}_1 (E_v \rho^0) = \mathcal{C}_m \cdots \mathcal{C}_1 \rho = \rho \]
  where in the last equality we have used that $\rho$ is acyclic.
\end{proof}

Note that the proof shows the following: if $\rho, \rho'$ are acyclic
infinitive stack configurations and $E_v \rho'=\rho$, then the unique
path in $\rho_0$ from $v$ to the sink is the \textbf{loop-erasure} of
the path taken by rotor-router walk started at~$v$ with initial
configuration~$\rho'$.

\section{Conjectures and Open Problems} \label{conc}

In this section we discuss some natural questions about chip-firing
and rotor-routing that remain unanswered.

Fey-den Boer and Redig \cite{FR} consider \textbf{aggregation} in the
sandpile model on~$\Z^d$.  In their setup, the underlying graph for
the chip-firing game is the infinite undirected $d$-dimensional cubic
lattice~$\Z^d$.  Start with each site containing $h \leq 2d-2$ chips.
Here $h$ may be even be taken negative, corresponding to starting with
a ``hole'' of depth $H=-h$ at each lattice site; that is, each site
absorbs the first $H$ chips it receives, and thereafter fires every
time it receives an additional four chips.  Now add $n$ chips to the
origin and stabilize.  Denote by $S_{n,H}$ the set of sites in~$\Z^d$
which fired in the process of stabilization.

\begin{thm}[\cite{FR}] \label{RF}
  Let~$\mathcal{C}(r)$ denote the cube of side length~$2r+1$ centered at the
  origin in $\Z^d$.  For each $n$ there exists an integer
  $r_n$ such that $S_{n,2-2d}=\mathcal{C}(r_n)$.
\end{thm}

In two dimensions, Theorem~\ref{RF} states that $S_{n,-2}$ is a
square.  Simulations indicate that for general $H\geq -2$,
the limiting shape of $S_{n,H}$ in $\Z^2$ may be a polygon with
$4H+12$ sides.

\begin{figure}[htbp]
\begin{center}
\begin{tabular}{ccc}
 \includegraphics[width=.3\textwidth]{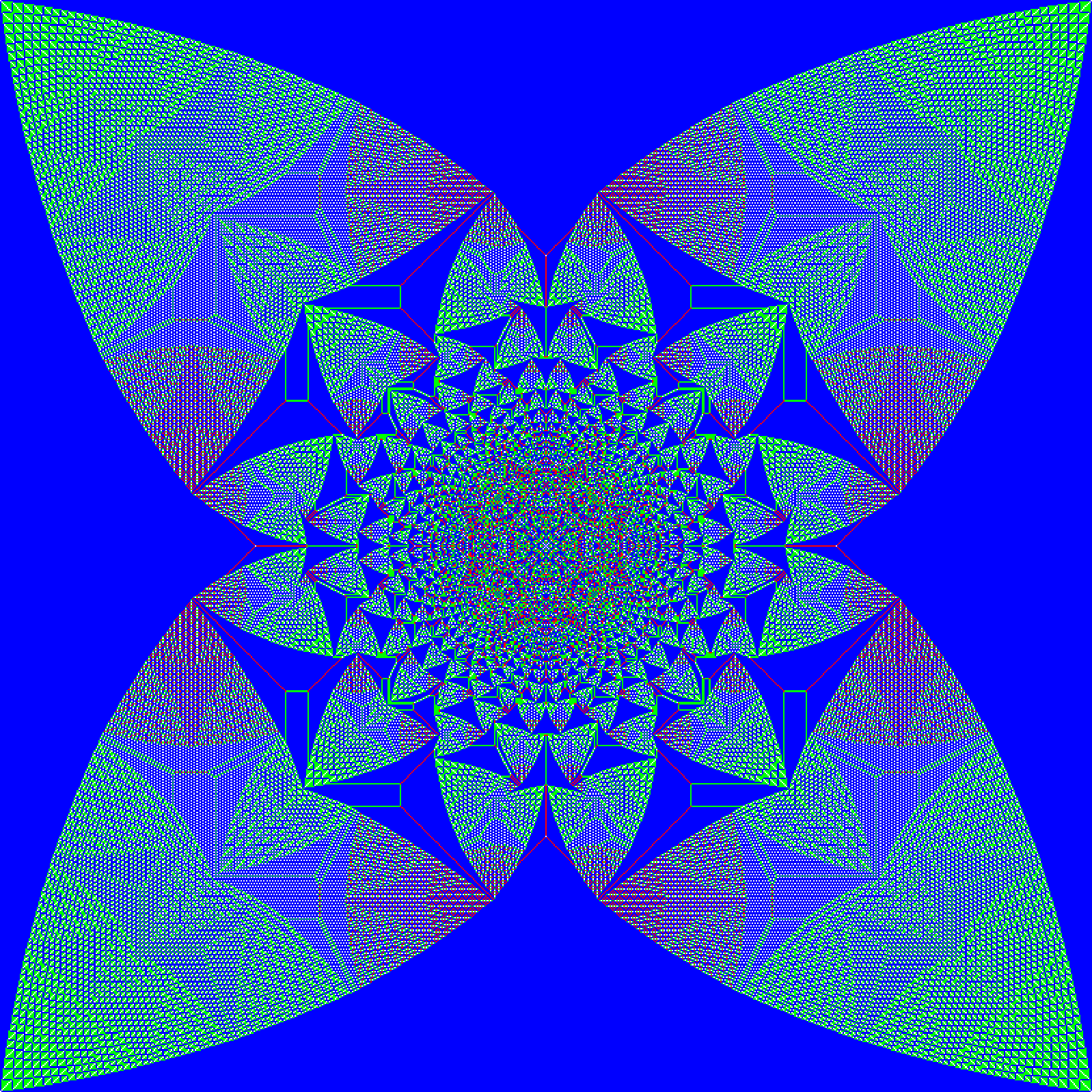}
&\includegraphics[width=.3\textwidth]{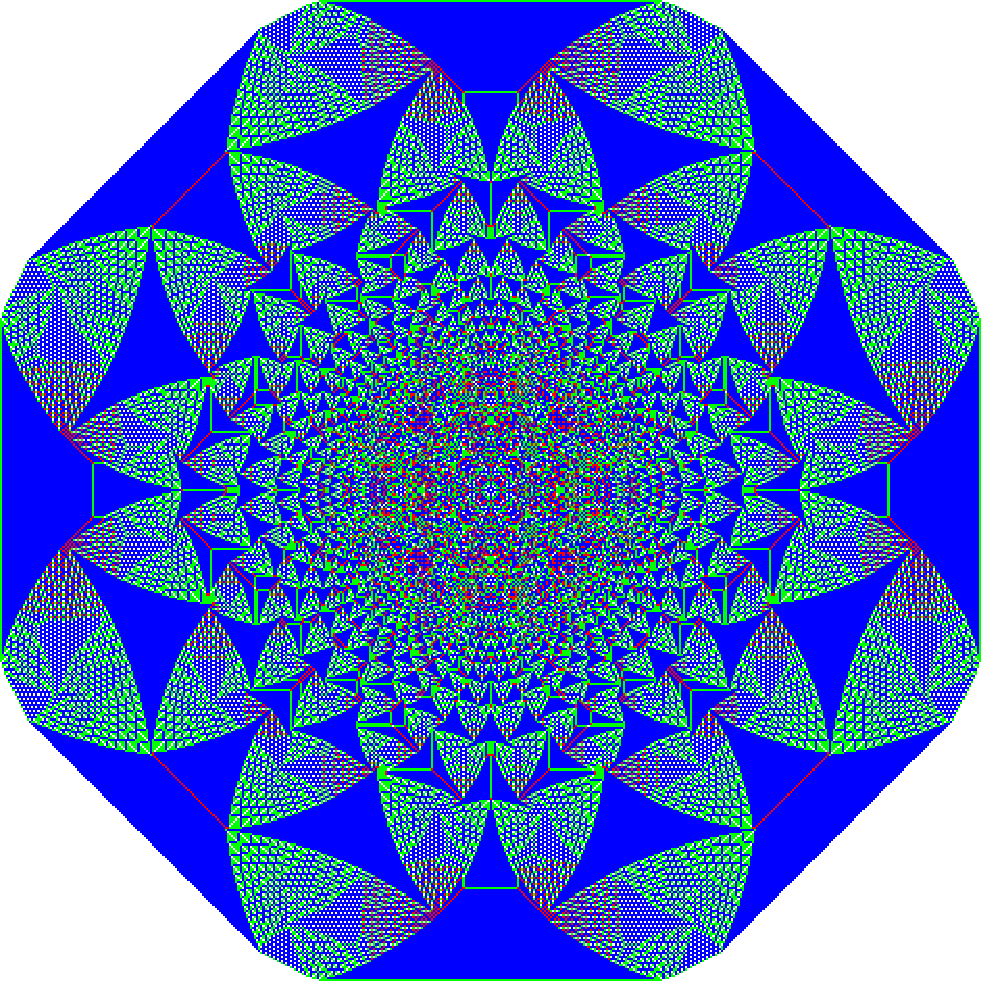}
&\includegraphics[width=.3\textwidth]{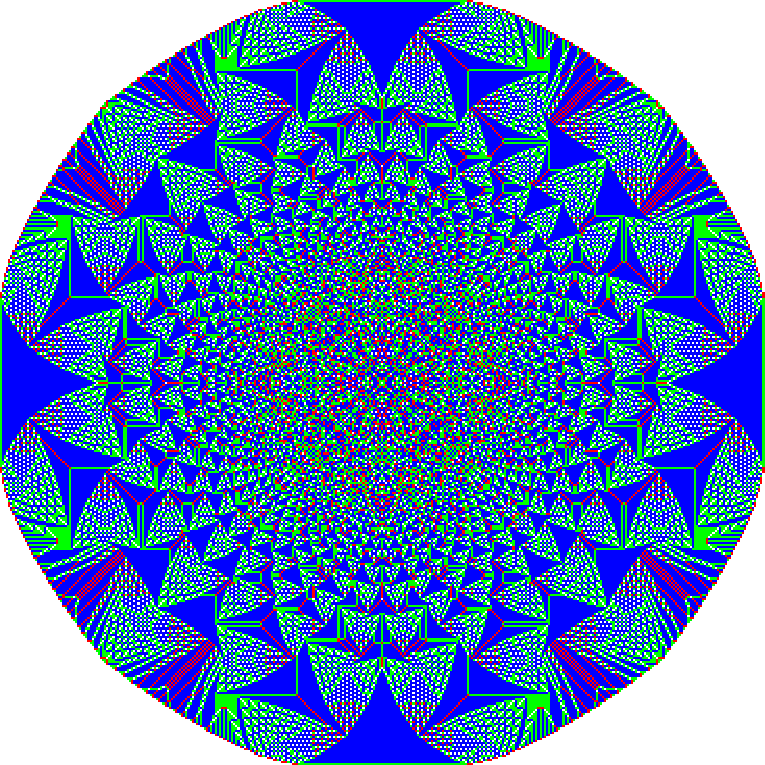} \\
$H=-2$ &$H=-1$ &$H=0$
\end{tabular}
\end{center}
\caption{Sandpile aggregates of $250,000$ chips in $\Z^2$
at hole depths $H=-2$ (left), $H=-1$ (center), and $H=0$.}
\label{polygons}
\end{figure}

\begin{q} \label{gon}
  In $\Z^2$, is the limiting shape of $S_{n,H}$ as $n \to \infty$ a
  regular $(4H+12)$-gon?  Simulations indicate a regular $(4H+12)$-gon
  with some ``rounding'' at the corners; it remains unclear if the
  rounded portions of the boundary become negligible in the limit.
  Even if the limiting shape is not a polygon, it would still be very
  interesting to establish the weaker statement that it has the
  dihedral symmetry $D_{4H+12}$.
\end{q}

The square, octagon and dodecagon corresponding to the cases
$H=-2,-1,0$ are illustrated in Figure~\ref{polygons}.  Regarding
Question~\ref{gon}, we should note that even the existence of a
limiting shape for $S_{n,H}$ has not been proved in the case $H>-2$.
On the other hand, as $H \to \infty$ the limiting shape is a ball in
all dimensions, as shown by Fey and Redig \cite{FR} and strengthened
in \cite{LP2}.  In the theorem below, $\omega_d$ denotes the volume of
the unit ball in $\R^d$, and $B_r$ denotes the discrete ball
        \[ B_r = \{x \in \Z^d \,|\, x_1^2 + \ldots + x_d^2 < r^2\}. \]

\begin{thm} \cite{LP2}.
  Fix an integer $H \geq 2-2d$.  Let $S_{n,H}$ be the set of sites in
  $\Z^d$ that fire in the process of stabilizing $n$ particles at the
  origin, if every lattice site begins with a hole of depth $H$.
  Write $n = \omega_d r^d$.  Then
  \[ B_{c_1r - c_2} \subset S_{n,H} \]
  where
  \[ c_1 = (2d-1+H)^{-1/d} \]
  and $c_2$ is a constant depending only on $d$.  Moreover if $H \geq
  1-d$, then for any $\epsilon>0$ we have
  \[ S_{n,H} \subset B_{c'_1r + c'_2} \]
  where
  \[ c'_1 = (d-\epsilon+H)^{-1/d} \]
  and $c'_2$ is independent of $n$ but may depend on $d$, $H$ and $\epsilon$.
\end{thm}

In particular, note that the ratio $c_1/c'_1 \uparrow 1$ as 
$H \uparrow \infty$.

For many classes of graphs, the identity element of the sandpile
group has remarkable properties that are not well understood.  Let
$I_n$ be the identity element of the $n\times n$ grid graph $G_n$
with wired boundary; the states $I_n$ for four different values of
$n$ are pictured in Figure~\ref{square-ident}.  Comparing the
pictures of $I_n$ for different values of $n$, one is struck by
their extreme similarity.  In particular, we conjecture that as $n
\to \infty$ the pictures converge in the following sense to a
limiting picture on the unit square $[0,1]\times [0,1]$.

\begin{conj}
  Let $a_n$ be a sequence of integers such that $a_n \uparrow \infty$
  and $\frac{a_n}{n} \downarrow 0$.  For $x \in [0,1]\times [0,1]$ let
  \[ f_n(x) =  \frac{1}{a_n^2} \sum_{\substack{y \in G_n \\ ||y-nx|| \leq a_n}} I_n(y). \]
  There is a sequence $a_n$ and a function $f : [0,1] \times [0,1] \to
  \R_{\geq 0}$ which is locally constant almost everywhere, such that $f_n \to
  f$ at all continuity points of $f$.
\end{conj}

Most intriguing is the apparent fractal structure in the
conjectural~$f$.  Recent progress has been made toward understanding
the fractal structure of the identity element of a certain orientation
of the square grid; see \cite{CPS}.

By Lemma~\ref{Claim 6}, the recurrent orbits of the rotor-router
operation on an Eulerian digraph are extremely short: although the
number of unicycles is typically exponential in the number of
vertices, the orbits are all of size equal to the number of edges.
One would expect that such short orbits are not the norm for
general digraphs.

\begin{q}
  Does there exist an infinite family of non-Eulerian strongly
  connected digraphs $G_n$, such that for each $n$, all the unicycles
  of $G_n$ lie in a single orbit of the rotor-router operation?
\end{q}

Another question stemming from Lemma~\ref{Claim 6} is the following.
Fix two edges $e_0$ and $e_1$ of a digraph $G$.  Starting from a unicycle
on $G$, record a $0$ each time the chip traverses the edge $e_0$, and
record a $1$ each time it traverses $e_1$.  If $G$ is Eulerian, then
Lemma~\ref{Claim 6} implies that the resulting sequence will simply
alternate $0,1,0,1,\ldots$.  For a general digraph, the sequence is
periodic, since the initial unicycle must recur; what can be said
about the period?

Lastly, the articles \cite{PDDK} and \cite{PPS} contain several
conjectures that are supported by both credible heuristics and
computer experiments, but that have not been rigorously proved.  For
instance, it appears that, with random initial rotor orientations, the
set of sites visited by a rotor-router walk of length $n$ in the plane
typically has diameter on the order of $n^{1/3}$ \cite{PDDK} (compare
this with the corresponding growth rate for random walk in the plane,
which is $n^{1/2}$).

\section*{Acknowledgment}

We thank Darij Grinberg for pointing out a number of typos in the 2008 version, which are corrected in this version.

\bibliographystyle{halpha}
\bibliography{sand}
%\newpage

\end{document}